\documentclass[12pt,draftcls,onecolumn]{IEEEtran}

\usepackage{amsmath,amssymb,amsfonts}
\usepackage{graphicx}
\usepackage{xcolor}
\usepackage{algorithm}
\usepackage{algorithmic}

\newtheorem{definition}{Definition}
\newtheorem{assumption}{Assumption}
\newtheorem{lemma}{Lemma}
\newtheorem{theorem}{Theorem}
\newtheorem{remark}{Remark}

\newcommand{\R}{\mathbb{R}}
\newcommand{\X}{\mathcal{X}}
\newcommand{\Y}{\mathcal{Y}}
\newcommand{\B}{\mathbb{B}}
\newcommand{\RegF}{\operatorname{G-Regret}^F_T}
\newcommand{\argmin}{\mathop{\mathrm{arg\,min}}}
\newcommand{\argmax}{\mathop{\mathrm{arg\,max}}}
\newcommand{\dist}{\mathrm{dist}}

\begin{document}

\title{Forgetting-Factor Regret for Online Zero-Sum Games}

\author{Yuhang Liu, Zi'ang Yan, Wenjun Mei and Wenxiao Zhao$^{*}$%
\thanks{Yuhang Liu and Zi'ang Yan contributed equally to this work. Yuhang Liu and Wenjun Mei are with the Department of Control and Systems Engineering, Peking University, Beijing 100871, China; Zi'ang Yan and Wenxiao Zhao are with the Key Laboratory of Systems and Control, Academy of Mathematics and Systems Science, Chinese Academy of Sciences, Beijing 100190, China, and also with the School of Mathematical Sciences, University of Chinese Academy of Sciences, Beijing 100049, China (e-mail: yuhang.liu@pku.edu.cn; ftqeyza@163.com; mei@pku.edu.cn; wxzhao@amss.ac.cn). $^{*}$Corresponding author. }}

\maketitle

\begin{abstract}
This paper studies dynamic equilibrium tracking in online two-player zero-sum games with time-varying convex-concave payoff functions. Existing regret metrics for online saddle-point problems usually aggregate historical payoffs with uniform weights, and hence may fail to characterize the real-time tracking performance with respect to the current Nash equilibrium (NE). To address this issue, we introduce a zero-sum game regret function with a forgetting factor, which assigns exponentially decaying weights to past saddle gaps and emphasizes recent performance. This metric directly links regret minimization to the tracking of time-varying NEs. Within this framework, we investigate three online algorithms under different computational and information settings. For first-order feedback, we analyze projected gradient descent-ascent and design a projection-free online Frank-Wolfe method to reduce the computational cost of projections. For zeroth-order feedback, we develop a deterministic finite-difference method that only uses function-value queries. For all three algorithms, we establish forgetting-factor regret bounds that explicitly characterize the effects of NE variation, payoff variation, and gradient-estimation error. We further provide sufficient conditions under which the proposed regret converges to zero, thereby certifying asymptotic tracking of time-varying NEs. The numerical example validates the theoretical results and illustrates  the tracking advantage of the proposed regret metric.

\end{abstract}

\begin{IEEEkeywords}
Online zero-sum games, forgetting-factor regret, Nash equilibrium tracking, gradient-free algorithms, projection-free algorithms.
\end{IEEEkeywords}

\section{Introduction}
\label{sec:introduction}
\IEEEPARstart{G}{ame} theory aims to develop mathematical models for scenarios involving multi-player competition or cooperation, which has broad applications across fields such as economics, mathematics, biology, computer science, etc  \cite{nash1950equilibrium,von2007theory}. While optimization theory \cite{2003Zinkevich,2007Hazan} focuses on single-agent decision-making to optimize a specific objective function, game theory addresses multi-agent decision-making, aiming to find a balanced strategy combination that satisfies the mutual constraints among players.

Two-player zero-sum games \cite{1928Neumann,2004Morgenstern,2009Ben} constitute one of the most fundamental models in game theory, which have been widely applied in artificial intelligence problems, such as boosting~\cite{1996Freund}, generative adversarial networks~\cite{2014Good}, and poker games~\cite{2015Bow}. In a two-player zero-sum game, one player minimizes the payoff function $f(x,y):\X\times\Y\rightarrow \R$, while the other player maximizes it. A point pair $(x^*,y^*)$ is called a Nash equilibrium (NE), or equivalently a saddle point, if
\begin{align}
f(x^*,y)\leq f(x^*,y^*)\leq f(x,y^*), \quad \forall (x,y)\in \X\times \Y. \label{eq:static_ne}
\end{align}
Kakutani's fixed-point theorem guarantees the existence of an NE for compact convex decision sets and continuous convex-concave payoff functions~\cite{Kakutani1941}. Many algorithms have been developed for saddle-point computation, including gradient descent-ascent algorithms~\cite{Nedic2009Saddle,daskalakis2018limit,schafer2019competitive,adolphs2019local}, primal-dual algorithms~\cite{hamedani2021primal}, extra-gradient and optimistic gradient methods~\cite{1995Tseng,2019Liang,2020Mokhtari,mokhtari2020convergence,yoon2021accelerated,2021weilinear}, gradient-free algorithms~\cite{liu2020min,wang2023zeroth}, Frank-Wolfe methods~\cite{gidel2017frank,abernethy2017frank,abernethy2018faster}, and nonconvex game algorithms~\cite{hazan2017efficient,lin2020gradient}.

Recently, online two-player zero-sum game, also known as online saddle-point problem, with time-varying convex-concave payoff functions has received increasing attention. This setting extends both online convex optimization~\cite{Shalevshwartz2012Online,2004Online,2012projection} and static zero-sum games. 
In such a scenario, the Nash equilibrium (saddle point)  seeking algorithms need to be performed for a sequence of convex-concave payoff functions $\{f_t(\cdot,\cdot):\X\times \Y\rightarrow \mathbb{R},t=0,\dots, T\}$. In each round $t$, the current strategy pair $(x_t,y_t)$ is evaluated under $f_t$, and the available information up to time $t$ is then used to generate $(x_{t+1},y_{t+1})$. 
Motivated by static regret \cite{Shalevshwartz2012Online,2003Zinkevich,2019Introduction} and dynamic regret \cite{2018Dynamic,7963457} in online convex optimization, various regret functions have been proposed for online two-player zero-sum
games with time-varying payoff functions. 
Literature \cite{2018Adrian} introduces the individual-regret
\begin{align}
& \operatorname{Ind-Regret}_{x}(T)=\sum_{t=1}^T f_t(x_t, y_t)-\min _{x \in \X} \sum_{t=1}^T f_t(x, y_t), \label{eq:ind_x}\\
& \operatorname{Ind-Regret}_{y}(T)=\max _{y \in \Y} \sum_{t=1}^T f_t(x_t, y)-\sum_{t=1}^T f_t(x_t, y_t), \label{eq:ind_y}
\end{align}
and establishes the saddle-point regret
\begin{align}
\operatorname{SP-Regret}_T
=\left|\sum_{t=1}^T f_t(x_t, y_t)-\min _{x \in \X} \max _{y \in \Y} \sum_{t=1}^T f_t(x, y)\right|. \label{eq:sp_reg}
\end{align}
It further proposes a saddle-point Follow-The-Leader algorithm and establishes sublinear saddle-point regret bounds. Building on this framework, \cite{2-xu2019online} introduces regularization terms into the saddle-point Follow-The-Leader algorithm to stabilize the estimates and establishes the sub-linear saddle-point regret bound. Additionally, \cite{28-cardoso2019competing} studies online bilinear games, where the payoff functions are time-varying matrices, and derives the corresponding saddle-point regret. 

Note that individual regret and saddle-point regret do not account for the influence of time. Therefore, algorithms with sublinear regrets do not necessarily guarantee satisfactory tracking performance for time-varying Nash equilibria.
{\em Example:} 
Let $T$ be even, $\X=\Y=[-2,2]$, and
\begin{equation}
f_t(x,y)=
\begin{cases}
((x-1)^2+5)(5-y^2), & t=0,\ldots,\frac{T}{2}-1,\\
((x+1)^2+5)(5-y^2), & t=\frac{T}{2},\ldots,T.
\end{cases}
\end{equation}

A direct calculation shows that the time-varying Nash equilibria are given by 
\begin{equation}\label{SFTL}
\begin{aligned}
(x_{t}^*,y_t^*)=\left\{\begin{array}{ll}
(1,0), & t=0, \dots, \frac{T}{2}-1, \\
(-1,0), & t=\frac{T}{2}, \dots, T.
\end{array}\right.
\end{aligned}
\end{equation}
By using the saddle-point Follow-The-Leader algorithm (see, e.g., \cite{2018Adrian}):
\begin{align*}
\left(x_{t+1},y_{t+1}\right)=\arg\min_{x\in \X}\max_{y\in \Y}\sum_{s=0}^t f_s(x, y),
\end{align*}
or by using Follow-The-Leader algorithm for online convex optimization (see, e.g., \cite{Shalevshwartz2012Online}, at time $t$, define the loss function by $f_t(x,y_t)$ with respect to variable $x$ and by $f_t(x_t,y)$ with respect to variable $y$), 
\begin{align*}
x_{t+1}&=\arg\min_{x\in \X}\sum_{s=0}^t f_s(x,y_s),\\
y_{t+1}&=\arg\max_{y\in \Y}\sum_{s=0}^t f_s(x_s,y).
\end{align*}
it follows that both algorithms generate 
\begin{equation}
(x_t,y_t)=
\begin{cases}
(1,0), & t=1,\ldots,\frac{T}{2},\\
(T/t-1,0), & t=\frac{T}{2}+1,\ldots,T,
\end{cases}
\end{equation}
The saddle-point regret bound of saddle-point Follow-The-Leader algorithm is $ O(\log T) $  (see, Theorem 1 in \cite{2018Adrian}). The individual regret $\operatorname{Ind-Regret}_x(T)$ bound of Follow-The-Leader algorithm is $ O(\log T) $ (see, Corollary 2.2 in \cite{Shalevshwartz2012Online}), 
 and it can be calculated directly that $\operatorname{Ind-Regret}_y(T)=0$. 
On the other hand, it is direct to check that the estimate $x_T$ equals $0$, rather than the Nash equilibrium $-1$.

Based on $\operatorname{SP-Regret}_T$ defined in (\ref{eq:sp_reg}), some other regrets are given in the literature for evaluating the performance of algorithms. 
\cite{3-roy2019online,29-zhang2022no} proposed the dynamic Nash equilibrium regret $\operatorname{NE-Regret}$ defined by
\begin{align}
\operatorname{NE-Regret}_T
=\left|\sum_{t=1}^T f_t(x_t, y_t)-\sum_{t=1}^T f_t(x_t^*, y_t^*)\right|, \label{eq:ne_reg}
\end{align}
where $(x_t^*,y_t^*)$ is a Nash equilibrium of payoff function $f_t(x,y)$ for $t\geq 1$, i.e., 
\begin{align}
f_t(x^*_t,y_t)\leq f_t(x^*_t,y^*_t)\leq f_t(x_t,y^*_t).\label{NE}
\end{align}
\cite{29-zhang2022no} introduces $\operatorname{NE-Regret}$ for online zero-sum games with bilinear payoff functions and gives the regret bound corresponding to a NE variation $U_T\triangleq \sum_{t=1}^T\left(\left\|x_{t}^*-x_{t-1}^*\right\|_1+\left\|y_{t}^*-y_{t-1}^*\right\|_1\right)$. 
\cite{3-roy2019online} investigates the general convex-concave payoff functions and obtains the regret bound corresponding to a NE variation  $V_T\triangleq \sum_{t=1}^T\left(\left\|x_{t}^*-x_{t-1}^*\right\|^2+\left\|y_{t}^*-y_{t-1}^*\right\|^2\right)$. 
\cite{29-zhang2022no,27-meng2023online} introduce the duality gap, which can be summarized by 
\begin{align}
\operatorname{Dual-Gap}_T(u_t, v_t)=
\sum_{t=1}^T f_t(x_t, v_t)-\sum_{t=1}^T f_t(u_t, y_t), \label{eq:dual_gap}
\end{align}
where $\{(u_t,v_t)\in \X\times \Y\}_{1\leq t\leq T}$ is an arbitrary sequence. \cite{27-meng2023online} provides a duality gap bound corresponding to the variation $ P_T \triangleq \sum_{t=1}^T \left( \|u_t - u_{t-1}\| + \|v_t - v_{t-1}\| \right) $. It's also mentioned that for a saddle point game, it suffices to choose $(u_t, v_t) = (x_t^*, y_t^*)$,  which then allows the derivation of a deterministic $\operatorname{Dual-Gap}_T(x_t^*, y_t^*)$ bound accounting for a NE variation.


In the regret notions above, the time-varying payoff functions $\{f_t(\cdot,\cdot),~t=1,\dots,T\}$ are aggregated with equal weights. In real-time decision systems, however, terminal or recent performance is usually more informative than historical average performance. For tracking a moving NE, the desired terminal property is that the two one-sided saddle gaps vanish, namely $f_T(x_T,y_T^*)-f_T(x_T^*,y_T^*)\to0$ and $f_T(x_T^*,y_T^*)-f_T(x_T^*,y_T)\to0$, or equivalently $f_T(x_T,y_T^*)-f_T(x_T^*,y_T)\to0$. 

However, neither $\operatorname{NE-Regret}_T$ in \eqref{eq:ne_reg} nor
$\operatorname{Dual-Gap}_T(x_t^*,y_t^*)$ in \eqref{eq:dual_gap} can guarantee
$f_t(x_t,y_t^*)-f_t(x_t^*,y_t)\to0$. The former may vanish merely because the played strategies attain the same payoff value as the NE without approaching the NE, whereas the latter may be sublinear even when large instantaneous saddle gaps occur infinitely often.
For $\operatorname{NE-Regret}_T$, we consider the above example and let $(x_t,y_t)=(2,\sqrt{5/6})$ for $t=0,\dots,T/2-1$ and  $(x_t,y_t)=(0,\sqrt{5/6})$ for $t=T/2,\dots,T$. It is direct to check $f_t\left(x_{t}, y_{t}\right)=f_t(x_t^*, y_t^*)$ for $t=1,\dots,T$ and thus $\operatorname{NE-Regret}_T=0$. However, the estimates $\{(x_t,y_t)\}_{t\ge 0}$ do not converge to the NE $(-1, 0)$, and $f_t(x_t, y_t^*) - f_t(x_t^*, y_t)$ does not converge to zero either. 
For the duality gap $\operatorname{Dual-Gap}_T(x_t^*, y_t^*)$, define a sequence $\{a_t\}_{t\geq 0}$ by $a_t=1$ if $t=2^m$, $m\geq1$ and $a_t=0$ otherwise. If $f_t(x_t, y_t^*) - f_t(x_t^*, y_t)=a_t$ for $t\geq 0$, then $\lim_{T\rightarrow \infty}\sum_{t=1}^T a_t/T\leq \lim_{T\rightarrow \infty}\log_2 T/T=0$, while $f_t(x_t, y_t^*) - f_t(x_t^*, y_t)$ does not converge to zero. 

Motivated by the forgetting-factor regret for online convex optimization~\cite{liu2023forgetting}, we introduce exponential weights into the stage-wise saddle gaps and propose the following online zero-sum game regret with a forgetting factor:
\begin{align}
\operatorname{G-Regret}^F_T
\triangleq \sum_{t=1}^T \rho^{T-t}\left(f_t(x_t, y_t^*)-f_t(x_t^*, y_t)\right), \label{eq:intro_gregret}
\end{align}
where $\rho\in(0,1)$ is the forgetting factor. {The exponent $T-t$ makes the contribution of an old game decay as its distance from the terminal time increases, while the current terminal gap has weight one. The term inside the sum is the two-sided saddle gap against the current NE, so the proposed metric directly evaluates whether the online decisions track the moving equilibrium rather than merely achieving a good historical average.}

Under the framework of the forgetting-factor game regret, this paper addresses two questions for online two-player zero-sum games: 
\begin{itemize}
\item[(i)] how large $\operatorname{G-Regret}^F_T$ is for concrete algorithms, and
\item[(ii)] under what verifiable variation conditions $\operatorname{G-Regret}^F_T$ tends to zero. 
\end{itemize}
We answer these questions for projected online gradient descent-ascent, projection-free Frank-Wolfe updates, and a deterministic-difference zeroth-order method. 
If $\operatorname{G-Regret}^F_T=o(1)$, then by nonnegativity and the NE condition,
\begin{align}
0\le f_T(x_T,y_T^*)-f_T(x_T^*,y_T^*)\le \operatorname{G-Regret}^F_T=o(1),\\
0\le f_T(x_T^*,y_T^*)-f_T(x_T^*,y_T)\le \operatorname{G-Regret}^F_T=o(1),
\end{align}
which establishes terminal tracking in payoff value. {Under strong convexity-strong concavity, this further implies convergence of $(x_T,y_T)$ to the time-varying NE in squared distance.}

The contributions of the paper are as follows.
\begin{itemize}
\item {We introduce a forgetting-factor regret for online zero-sum games and show that its convergence to zero directly certifies real-time tracking of the time-varying NE. This provides a performance measure that is stronger than historical average guarantees for tracking purposes.}
\item We analyze the performance of {projected online gradient descent-ascent} for online two-player zero-sum {games} under $\operatorname{G-Regret}^F_T$. {A sharp contraction factor is derived under strong convexity-strong concavity and cross-gradient Lipschitz conditions, together with sufficient conditions on the NE variation for $\operatorname{G-Regret}^F_T=o(1)$.}
\item {We design and analyze a projection-free online Frank-Wolfe method for games. The regret bound reveals how cross-gradient coupling enters both the contraction condition and the admissible adaptive step size, while payoff variation and NE variation appear as explicit disturbance terms.}
\item {We further develop a deterministic-difference zeroth-order algorithm for the setting where gradients are unavailable. Its regret bound separates the effects of equilibrium variation and gradient-estimation accuracy, and shows that vanishing perturbation radii preserve tracking.}
\end{itemize}

The rest of the paper is organized as follows. {Section~II formulates the online zero-sum game and the forgetting-factor regret. Section~III studies first-order algorithms, including projected gradient descent-ascent and Frank-Wolfe updates. Section~IV gives the zeroth-order extension. Section~V presents simulations, Section~VI concludes the paper, and the Appendix contains auxiliary lemmas and proofs.}

{\em Notation.} Denote by $\|\cdot\|$ the Euclidean norm and by $\langle\cdot,\cdot\rangle$ the Euclidean inner product. For a closed convex set $\X$, $\mathcal{P}_{\X}(z)\triangleq \argmin_{x\in\X}\|x-z\|$ denotes the Euclidean projection. Let $e_x^k$ and $e_y^\ell$ be the canonical basis vectors in $\R^{n_x}$ and $\R^{n_y}$, respectively. For a differentiable payoff $f(x,y)$, $\nabla_x f$ and $\nabla_y f$ denote its partial gradients.

\section{Problem Formulation}
\label{sec:problem}

Consider online two-player zero-sum games with payoff functions
\[
f_t:\X\times\Y\rightarrow\R,\qquad t=0,1,\ldots,
\]
where $\X\subset\R^{n_x}$ and $\Y\subset\R^{n_y}$ are feasible strategy sets. At time $t$, the minimizing player selects $x_t\in\X$ and the maximizing player selects $y_t\in\Y$. The pair $(x_t^*,y_t^*)$ is an NE of the $t$-th game if
\begin{align}
f_t(x_t^*,y)\le f_t(x_t^*,y_t^*)\le f_t(x,y_t^*),\quad \forall (x,y)\in\X\times\Y. \label{eq:online_ne}
\end{align}
\begin{definition}[Forgetting-factor game regret]
\label{def:ff_game_regret}
Let $\mathcal A$ be an online algorithm for the above zero-sum game, and denote by $\{(x_t,y_t)\}_{t=0}^{T}$ the strategy pairs generated by $\mathcal A$.
For a given $\rho\in(0,1)$, the forgetting-factor game regret of
$\mathcal A$ is defined by
\begin{align}
\RegF
\triangleq \sum_{t=1}^T \rho^{T-t}
\left(f_t(x_t,y_t^*)-f_t(x_t^*,y_t)\right). \label{eq:regret_def}
\end{align}
\end{definition}

In \eqref{eq:regret_def}, the coefficient $\rho^{T-t}$ decreases as the gap between $t$ and the terminal time $T$ increases. Hence recent games receive larger weights, and past games are progressively forgotten. This is consistent
with the tracking objective: when the environment varies, the quality of $(x_T,y_T)$ relative to the current NE is more relevant than the accumulated performance over obsolete games, i.e., 
\begin{align}
f_T(x_T,y_T^*)-f_T(x_T^*,y_T)\le \RegF. \label{eq:terminal_gap_bound}
\end{align}
Moreover, by \eqref{eq:online_ne}, each summand in \eqref{eq:regret_def} is nonnegative and 
\begin{align}
f_t(x_t,y_t^*)-f_t(x_t^*,y_t)= &f_t(x_t,y_t^*)-f_t(x_t^*,y_t^*) \nonumber\\
&+ f_t(x_t^*,y_t^*)-f_t(x_t^*,y_t). \label{eq:gap_decomp}
\end{align}

Therefore, if $\RegF=o(1)$, then the terminal two-sided saddle gap vanishes. Under strong convexity-strong concavity, this also implies that $(x_T,y_T)$ approaches the time-varying NE in squared distance.

We impose the following assumptions.
\begin{assumption}
\label{ass:set}
The sets $\X$ and $\Y$ are nonempty, compact, and convex. Each $f_t$ is continuous and convex-concave on $\X\times\Y$.
\end{assumption}

Since $\X$ and $\Y$ are compact, their diameters are finite. Throughout the paper, we use
\begin{align}
D_\X&\triangleq\sup_{x,x'\in\X}\|x-x'\|,\qquad
D_\Y\triangleq\sup_{y,y'\in\Y}\|y-y'\|. \label{eq:diameters}
\end{align}

\begin{assumption}
\label{ass:strong}
There exist $\mu_x,\mu_y>0$ such that each $f_t$ is $\mu_x$-strongly convex in $x$ for every fixed $y\in\Y$, and $\mu_y$-strongly concave in $y$ for every fixed $x\in\X$.
\end{assumption}

\begin{assumption}
\label{ass:smooth}
Each $f_t$ is differentiable, and there exist $L_x\ge\mu_x,L_y\ge\mu_y$ such that
\begin{align}
\|\nabla_x f_t(x,y)-\nabla_x f_t(x',y)\|&\le L_x\|x-x'\|,\label{eq:own_x}\\
\|\nabla_y f_t(x,y)-\nabla_y f_t(x,y')\|&\le L_y\|y-y'\| \label{eq:own_y}
\end{align}
for all $(x,y)\in \X\times\Y$ and $t\ge 0$.
\end{assumption}

By Lemma \ref{Lem:sigma} in  Appendix \ref{Appendix:A}, Assumption \ref{ass:smooth} is equivalent to that $f_t$ is $L_x~(L_y)-$smooth with respect to $x~(y)$ for any fixed $y~(x)$. 
\begin{assumption}
\label{ass:cross}
There exist $L_{xy},L_{yx}>0$ such that
\begin{align}
\|\nabla_x f_t(x,y)-\nabla_x f_t(x,y')\|&\le L_{xy}\|y-y'\|,\label{eq:cross_xy}\\
\|\nabla_y f_t(x,y)-\nabla_y f_t(x',y)\|&\le L_{yx}\|x-x'\| \label{eq:cross_yx}
\end{align}
for all $(x,y)\in \X\times\Y$ and $t\ge 0$.
\end{assumption}

\begin{remark}
Assumption~\ref{ass:cross} is peculiar to the game setting, as it controls the cross-player coupling in the gradients of the objective functions. In its absence, the influence of one player's strategy on the other player's gradient can be arbitrarily large, potentially overwhelming the gradient variation induced by the player's own strategy. Consequently, the gradient may fail to provide a meaningful local descent direction for that player. Such a condition is therefore crucial for establishing the validity and convergence of NE seeking algorithms.
\end{remark}

Before turning to the algorithmic analysis, we give two technical lemmas that are central to the subsequent development. The first lemma guarantees that the time-varying NE trajectory is well defined and unique under the
standing assumptions. The second lemma gives a refined gradient estimate for smooth strongly convex functions, which enables the sharp contraction analysis of the projected gradient descent-ascent update. These lemmas  constitute the basic analytical building blocks of our forgetting-factor regret bounds.
\begin{lemma}
\label{lem:exist_unique}
Under Assumptions~\ref{ass:set} and~\ref{ass:strong}, each $f_t$ admits a unique NE.
\end{lemma}

\begin{lemma}
\label{lem:refined}
Let function $\phi(\cdot)$ be differentiable, $\mu$-strongly convex, and $L$-smooth on a convex set $\X$ with nonempty interior. Then for any $u,~v\in \X$, it follows that 
\begin{align}
\|\nabla \phi(u)-\nabla \phi(v)-\frac{L+\mu}{2}(u-v)\|
\le \frac{L-\mu}{2}\|u-v\|. \label{eq:refined}
\end{align}
\end{lemma}
{\quad~~\it{Proof of Lemma \ref{lem:exist_unique}--\ref{lem:refined}: }}
See Appendix \ref{Appendix:B}.

Lemma~\ref{lem:exist_unique} allows us to treat the NE sequence as a
single-valued reference trajectory. For $t\ge1$, define
\begin{equation}
\begin{split}\label{eq:variation}
&V_t^x\triangleq \|x_t^*-x_{t-1}^*\|,\\
&V_t^y\triangleq \|y_t^*-y_{t-1}^*\|.
\end{split}
\end{equation}
These quantities measure the temporal variation of the equilibrium trajectory.
Since exact tracking is impossible in the online setting when the NE moves arbitrarily fast, they serve as disturbance terms in the subsequent tracking and regret bounds.

\section{First-Order Online Algorithms}
\label{sec:first_order}

\subsection{Projected Online Gradient Descent-Ascent}
\label{subsec:ogda}

The projected online gradient descent-ascent (OGDA) update is one of the most widely applied algorithms for online zero-sum games, which can be formulated by Algorithm \ref{alg:ogda}. 
\begin{algorithm}[htbp]
\caption{Projected Online Gradient Descent-Ascent Algorithm}
\label{alg:ogda}
\begin{algorithmic}[1]
\STATE {\bf  Initialization: } An initial estimate $(x_{0},y_0)$, a step size $\alpha>0$ and the maximal number $T$ of iterations.
\FOR{$t=0,\ldots,T$}
\STATE gradient calculation:
\begin{align*}
r_t^x&=\nabla_x f_t(x_t,y_t),\\
r_t^y&=\nabla_y f_t(x_t,y_t).
\end{align*}
\STATE projected gradient descent-ascent:
\begin{align}
x_{t+1}&=\mathcal{P}_{\X}(x_t-\alpha r_t^x),
\\
y_{t+1}&=\mathcal{P}_{\Y}(y_t+\alpha r_t^y). 
\end{align}
\ENDFOR
\end{algorithmic}
\end{algorithm}

Next we analyze the forgetting-factor regret  for Algorithm \ref{alg:ogda}. We first define a set of constants to be used
later on:
\begin{equation}
\begin{split}
\mu&\triangleq \min\{\mu_x,\mu_y\},\\
L&\triangleq \max\{L_x,L_y\},\\
L_\times&\triangleq \max\{L_{xy},L_{yx}\}, \label{eq:muL}
\end{split}
\end{equation}
and
\begin{align}
\beta(\alpha)
&\triangleq 1+\frac{1}{2}(L^2+\mu^2+2L_\times^2)\alpha^2-(L+\mu)\alpha \nonumber\\
&\quad +\frac{1}{2}(L-\mu)\alpha
\sqrt{(L\alpha+\mu\alpha-2)^2+4L_\times^2\alpha^2}, \label{eq:beta_gda}
\end{align}
where $\mu,L,L_\times$ are defined in \eqref{eq:muL}, and $\alpha$ is the step-size in Algorithm~\ref{alg:ogda}. Since  $L\ge \mu$, it can be calculated that if
\begin{align}
\alpha\in\left(0,\frac{2\mu}{\mu L+L_\times^2}\right) \label{eq:alpha_range}
\end{align}
holds, then 
\begin{align}
\beta(\alpha)\in[0,1). \label{eq:beta}
\end{align}
\begin{remark}
The contraction factor in \eqref{eq:beta_gda} explicitly captures the stabilizing effect of strong convexity-strong concavity and the destabilizing effect of cross-gradients. The admissible range \eqref{eq:alpha_range} is the condition under which the former dominates the latter. The minimizer of $\beta(\alpha)$ is
\begin{align}
\alpha^*
=\frac{2\big((L+\mu)\sqrt{L_\times^2+L\mu}-L_\times(L-\mu)\big)}
{\big(4L_\times^2+(L+\mu)^2\big)\sqrt{L_\times^2+L\mu}}.
\end{align}
At this step size, the minimum contraction factor is
\begin{align}
\beta_{\min}
&=\frac{8L_\times(L^2-\mu^2)\sqrt{L\mu+L_\times^2}
+(L^2-\mu^2)^2}{\big((L+\mu)^2+4L_\times^2\big)^2}
\nonumber\\
&\quad+\frac{16L_\times^2(L\mu+L_\times^2)}
{\big((L+\mu)^2+4L_\times^2\big)^2}.
\end{align}
\end{remark}
\begin{theorem}
\label{thm:ogda}
Suppose Assumptions~\ref{ass:set}--\ref{ass:cross} hold and the NE trajectory lies in the interior of $\X\times\Y$. Let $\{(x_t,y_t)\}_{t\ge 0}$ be generated by
Algorithm~\ref{alg:ogda} with a step size $\alpha$ satisfying
\eqref{eq:alpha_range}. Then, for $t\ge0$,
\begin{align}
&\|x_t-x_t^*\|^2+\|y_t-y_t^*\|^2 \nonumber \le \beta^{t}\big(\|x_0-x_0^*\|^2+\|y_0-y_0^*\|^2\big) \nonumber\\
&~+
\sum_{i=1}^t\beta^{t-i}
\big(2D_\X V_i^x+2D_\Y V_i^y +(V_i^x)^2+(V_i^y)^2\big), \label{eq:ogda_compact_e}
\end{align}
where $\beta=\beta(\alpha)$ is given by \eqref{eq:beta_gda}, $D_\X,D_\Y$ are
given by \eqref{eq:diameters}, and $V_i^x,V_i^y$ are given by
\eqref{eq:variation}. The corresponding forgetting-factor regret satisfies
\begin{align}
\RegF
&\le \frac{L}{2}\sum_{t=1}^{T}\rho^{T-t}\beta^{t}
\big(\|x_0-x_0^*\|^2+\|y_0-y_0^*\|^2\big)
\nonumber\\
&\quad+\frac{L}{2}\sum_{i=1}^{T}\sum_{t=i}^{T}
\rho^{T-t}\beta^{t-i}
\big(2D_\X V_i^x+2D_\Y V_i^y \nonumber\\
&\quad+(V_i^x)^2+(V_i^y)^2\big). \label{eq:ogda_compact_regret}
\end{align}
In particular, if $V_t^x,V_t^y=o(1)$, then $\RegF=o(1)$.
\end{theorem}
{\quad~~\it{Proof: }}See Appendix \ref{Appendix:C}. 

\subsection{Projection-Free Online Frank-Wolfe}
\label{subsec:fw}
Algorithm \ref{alg:ogda} requires projecting the iterate onto the feasible set. When the feasible set is high-dimensional, computing the projection typically involves solving a quadratic programming subproblem, which can be computationally expensive. 

Motivated by projection-free methods for online optimization and saddle-point problems~\cite{hazan2012online,gidel2017frank}, we consider the Frank--Wolfe-type online update under gradient feedback, which is given by Algorithm \ref{alg:fw}. Our goal here is to analyze this update under the forgetting-factor game regret introduced above. 



For $t\ge1$, define
\begin{align}
F_t^{\sup}\triangleq \sup_{(x,y)\in\X\times\Y}|f_t(x,y)-f_{t-1}(x,y)|. \label{eq:Fsup}
\end{align}

In addition, for the case of the NE trajectory remains uniformly in the interior of
$\X\times\Y$, define
\begin{align}
\delta_\X&\triangleq \inf_{t\ge0}\dist(x_t^*,\partial\X)>0,\nonumber\\
\delta_\Y&\triangleq \inf_{t\ge0}\dist(y_t^*,\partial\Y)>0,\label{def:delta}
\end{align}
and
\begin{align}
C&\triangleq L_xD_\X^2+L_yD_\Y^2+(L_{xy}+L_{yx})D_\X D_\Y,\label{def:C}\\
\nu&\triangleq
1-\frac{\sqrt{2}\max\{D_\X L_{xy}/\sqrt{\mu_y},D_\Y L_{yx}/\sqrt{\mu_x}\}}
{\min\{\sqrt{\mu_y}\delta_\Y,\sqrt{\mu_x}\delta_\X\}}, \label{eq:nu}
\end{align}
where $D_\X,D_\Y$ are constants given by \eqref{eq:diameters}. 

\begin{algorithm}[htbp]
\caption{Projection-Free Online Frank-Wolfe Algorithm}
\label{alg:fw}
\begin{algorithmic}[1]
\STATE {\bf Initialization:} An initial estimate $(x_0,y_0)$, the maximal number $T$ of iterations and constants $\nu,~C$ given in \eqref{def:C}--\eqref{eq:nu}.
\FOR{$t=0,\ldots,T$}
\STATE gradient calculation:
\begin{align*}
r_t^x&=\nabla_x f_t(x_t,y_t),\\
r_t^y&=\nabla_y f_t(x_t,y_t).
\end{align*}
\STATE linear oracle calculation:
\begin{align*}
s_t^x&=\argmin_{x\in\X}\langle r_t^x,x\rangle,\\
s_t^y&=\argmax_{y\in\Y}\langle r_t^y,y\rangle.
\end{align*}
\STATE gap calculation:
\begin{align*}
g_t^x&=\langle r_t^x,x_t-s_t^x\rangle,\\
g_t^y&=-\langle r_t^y,y_t-s_t^y\rangle,\\
g_t&=g_t^x+g_t^y.
\end{align*}
\STATE Frank-Wolfe descent-ascent:
\begin{align*}
\eta_t&=\frac{\nu g_t}{C},\\
x_{t+1}&=(1-\eta_t)x_t+\eta_t s_t^x,\\
y_{t+1}&=(1-\eta_t)y_t+\eta_t s_t^y.
\end{align*}
\ENDFOR
\end{algorithmic}
\end{algorithm}

\begin{theorem}
\label{thm:fw}
Suppose Assumptions~\ref{ass:set}--\ref{ass:cross} hold, the NE trajectory
remains uniformly in the interior of $\X\times\Y$, and $\nu>0$. Let
$\{(x_t,y_t)\}_{t\ge0}$ be generated by Algorithm~\ref{alg:fw}. Then
$\eta_t\in[0,1]$ for all $t\ge0$. Moreover, for all $t\ge0$, 
\begin{align}
&f_t(x_t,y_t^*)-f_t(x_t^*,y_t) \nonumber\\
&\le \beta_{\rm FW}^{t}
\left(f_0(x_0,y_0^*)-f_0(x_0^*,y_0)\right) \nonumber\\
&\quad+\sum_{i=1}^{t}\beta_{\rm FW}^{t-i}
\Big(D_\X L_{yx}V_i^y+D_\Y L_{xy}V_i^x \nonumber\\
&\quad+\frac{L_y}{2}(V_i^y)^2+\frac{L_x}{2}(V_i^x)^2
+2F_i^{\sup}\Big), \label{eq:fw_gap_bound}
\end{align}
and 
\begin{align}
\RegF
&\le \sum_{t=1}^{T}\rho^{T-t}\beta_{\rm FW}^{t}
\left(f_0(x_0,y_0^*)-f_0(x_0^*,y_0)\right) \nonumber\\
&\quad+\sum_{i=1}^{T}\sum_{t=i}^{T}
\rho^{T-t}\beta_{\rm FW}^{t-i}
\Big(D_\X L_{yx}V_i^y+D_\Y L_{xy}V_i^x \nonumber\\
&\quad+\frac{L_y}{2}(V_i^y)^2+\frac{L_x}{2}(V_i^x)^2
+2F_i^{\sup}\Big), \label{eq:fw_regret}
\end{align}
where
\begin{align}
\beta_{\rm FW}
&\triangleq 1-\frac{\min\{\mu_y\delta_\Y^2,\mu_x\delta_\X^2\}\nu^2}{C}\in (0,1).
\end{align}
Here $V_i^x,V_i^y$ are given by \eqref{eq:variation}, and $F_i^{\sup}$ is
given by \eqref{eq:Fsup}. If $V_t^x,V_t^y,F_t^{\sup}=o(1)$, then
$\RegF=o(1)$.
\end{theorem}
{\quad~~\it{Proof: }}See Appendix \ref{Appendix:D}. 

\section{Zeroth-Order Online Algorithm}
\label{sec:zeroth}

We next consider the case where exact gradients are unavailable and only function values can be queried.
\begin{assumption}
\label{ass:zeroth}
There exist $\bar\delta_x,\bar\delta_y>0$ such that $f_t$ can be queried on $(\X+\bar\delta_x\B)\times(\Y+\bar\delta_y\B)$ for all $t$. 
Moreover, the differentiability and  smoothness bounds in Assumption~\ref{ass:smooth} hold on this expanded domain with the same constants. 
\end{assumption}

Let $\{c_t\}_{t\ge0}$ be a positive sequence. At time $t$, choose
\begin{align}
\delta_t=\min\left\{\frac{c_t}{n_xL_x},\frac{c_t}{n_yL_y},\bar\delta_x,\bar\delta_y\right\}. \label{eq:delta_choice}
\end{align}
For $k=1,\ldots,n_x$ and $\ell=1,\ldots,n_y$, define the deterministic central-difference estimators
\begin{align}
h_t^x[k]&=\frac{f_t(x_t+\delta_t e_x^k,y_t)-f_t(x_t-\delta_t e_x^k,y_t)}{2\delta_t},\\
h_t^y[\ell]&=\frac{f_t(x_t,y_t+\delta_t e_y^\ell)-f_t(x_t,y_t-\delta_t e_y^\ell)}{2\delta_t}. \label{eq:zero_estimator}
\end{align}

\begin{algorithm}[htbp]
\caption{Deterministic-Difference Zeroth-Order Descent-Ascent Algorithm}
\label{alg:zero}
\begin{algorithmic}[1]
\STATE {\bf Initialization:} An initial estimate $(x_0,y_0)$, a step size $\alpha>0$, a positive sequence $\{c_t\}$ and the maximal number $T$ of iterations.
\FOR{$t=0,\ldots,T$}
\STATE deterministic-difference:
\\$\quad$ for $k=1,\ldots, n_x$:
\begin{align*}
h_t^x[k]&=\frac{f_t(x_t+\delta_t e_x^k,y_t)-f_t(x_t-\delta_t e_x^k,y_t)}{2\delta_t}
\end{align*}
\\$\quad$ for $l=1,\ldots, n_y$:
\begin{align*}
h_t^y[\ell]&=\frac{f_t(x_t,y_t+\delta_t e_y^\ell)-f_t(x_t,y_t-\delta_t e_y^\ell)}{2\delta_t}.
\end{align*}
\STATE zeroth-order descent-ascent:
\begin{align*}
x_{t+1}&=\mathcal{P}_{\X}(x_t-\alpha h_t^x),\\
y_{t+1}&=\mathcal{P}_{\Y}(y_t+\alpha h_t^y).
\end{align*}
\ENDFOR
\end{algorithmic}
\end{algorithm}

\begin{theorem}
\label{thm:zero}
Suppose Assumptions~\ref{ass:set}--\ref{ass:cross} and~\ref{ass:zeroth} hold and the NE trajectory lies in the interior of $\X\times\Y$. Let $\{(x_t,y_t)\}_{t=0}^{T+1}$ be generated by
Algorithm~\ref{alg:zero} with a step size $\alpha$ satisfying
\eqref{eq:alpha_range}. Then, for $t\ge0$,
\begin{align}
&\|x_t-x_t^*\|^2+\|y_t-y_t^*\|^2 \nonumber\\
&\le \beta^{t}\big(\|x_0-x_0^*\|^2+\|y_0-y_0^*\|^2\big) \nonumber\\\nonumber
&\quad+\sum_{i=1}^{t}\beta^{t-i}
\Big(2D_\X(\alpha c_{i-1}+V_i^x)
+2D_\Y(\alpha c_{i-1}+V_i^y) \\
&\quad+2\alpha^2c_{i-1}^2+(V_i^x)^2 +(V_i^y)^2\Big), \label{eq:zero_compact_e}
\end{align}
where $\beta=\beta(\alpha)$ is given by \eqref{eq:beta_gda}, $D_\X,D_\Y$ are
given by \eqref{eq:diameters}, and $V_i^x,V_i^y$ are given by
\eqref{eq:variation}.
Moreover,
\begin{align}
&\RegF
\le \frac{L}{2}\sum_{t=1}^{T}\rho^{T-t}\beta^{t}
\big(\|x_0-x_0^*\|^2+\|y_0-y_0^*\|^2\big) \nonumber\\
&~+\frac{L}{2}\sum_{i=1}^{T}\sum_{t=i}^{T}
\rho^{T\!-\!t}\beta^{t\!-\!i} \Big(2D_\X(\alpha c_{i\!-\!1}\!+\!V_i^x)
\!+\!2D_\Y(\alpha c_{i\!-\!1}\!+\!V_i^y) \nonumber\\
&~+2\alpha^2c_{i-1}^2+(V_i^x)^2 +(V_i^y)^2\Big). \label{eq:zero_compact_regret}
\end{align}
In particular, if $c_t,V_t^x,V_t^y=o(1)$, then $\RegF=o(1)$.
\end{theorem}
{\quad~~\it{Proof: }}See Appendix \ref{Appendix:E}. 

\section{Simulation Studies}
\label{sec:simulation}
\setcounter{example}{0}
We 
consider a time-varying zero-sum game with payoff function
\begin{equation}
f_t(x,y)=
\begin{cases}
\frac{1}{2}(x-1)^2+\frac{1}{10}(x-1)y-\frac{1}{2}y^2,
& t=0,\ldots,\frac{T}{2}-1,\\
\frac{1}{2}(x+1)^2+\frac{1}{10}(x+1)y-\frac{1}{2}y^2,
& t=\frac{T}{2},\ldots,T.
\end{cases}
\end{equation}
Let $\X=\Y=[-10,10]$. It is straightforward to verify that 
the time-varying NE is
\begin{equation*}
(x_t^*,y_t^*)=
\begin{cases}
(1,0), & t=0,\ldots,\frac{T}{2}-1,\\
(-1,0), & t=\frac{T}{2},\ldots,T, 
\end{cases}
\end{equation*}
and $\mu=\mu_x=\mu_y=1$, $L=L_x=L_y=1$, and $L_{\times}=L_{xy}=L_{yx}=0.1$. 

We compare saddle-point Follow-The-Leader, Algorithm~\ref{alg:ogda},
Algorithm~\ref{alg:fw}, and Algorithm~\ref{alg:zero} from the initial point
$(0,-1)$. Choose $T=100$ and $\rho=0.5$. We use $\alpha=1$ for Algorithms~\ref{alg:ogda} and~\ref{alg:zero}, and set $c_t=1/(t+1)$ for Algorithm~\ref{alg:zero}, which can be verified to satisfy the step-size conditions in Theorem \ref{thm:ogda} and \ref{thm:zero}. The saddle-point Follow-The-Leader algorithm is given by
\begin{align}
(x_{t+1},y_{t+1})
=\argmin_{x\in\X}\max_{y\in\Y}\sum_{s=0}^{t}f_s(x,y),
\end{align}
and its trajectory in this example is
\begin{align}
(x_t,y_t)
=
\begin{cases}
(1,0), & t=1,\ldots,T/2,\\
(T/t-1,0), & t=T/2+1,\ldots,T.
\end{cases} \label{eq:sftl_path}
\end{align}

\begin{figure}[htbp]
\centering
\includegraphics[width=0.95\linewidth]{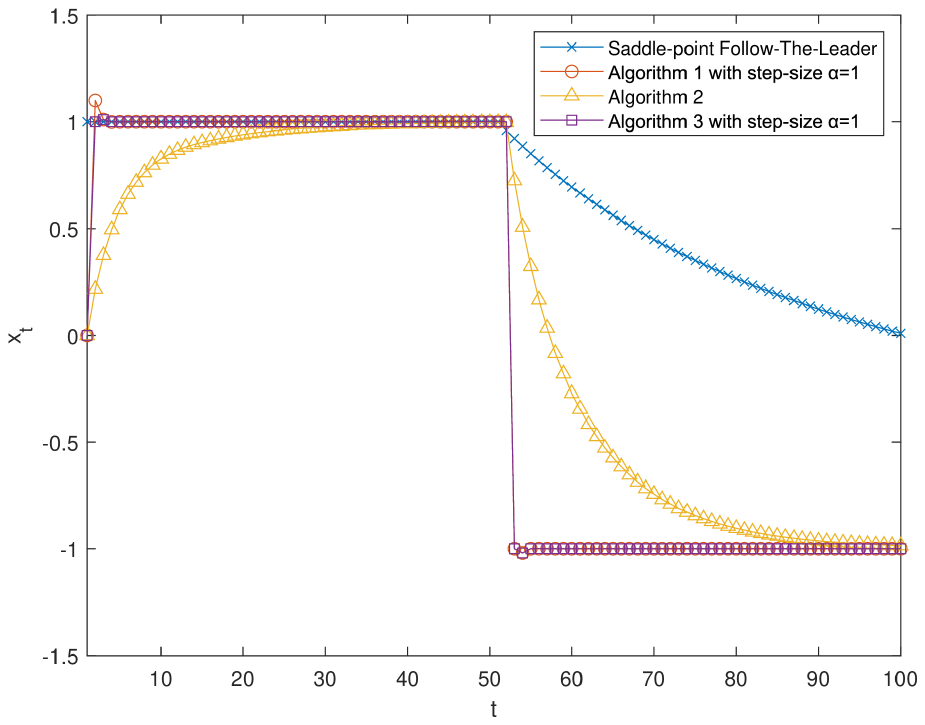}
\caption{Estimates of saddle-point Follow-The-Leader and Algorithms 1--3 with respect to $x$.}
\label{fig:x_first_order}
\end{figure}

\begin{figure}[htbp]
\centering
\includegraphics[width=0.95\linewidth]{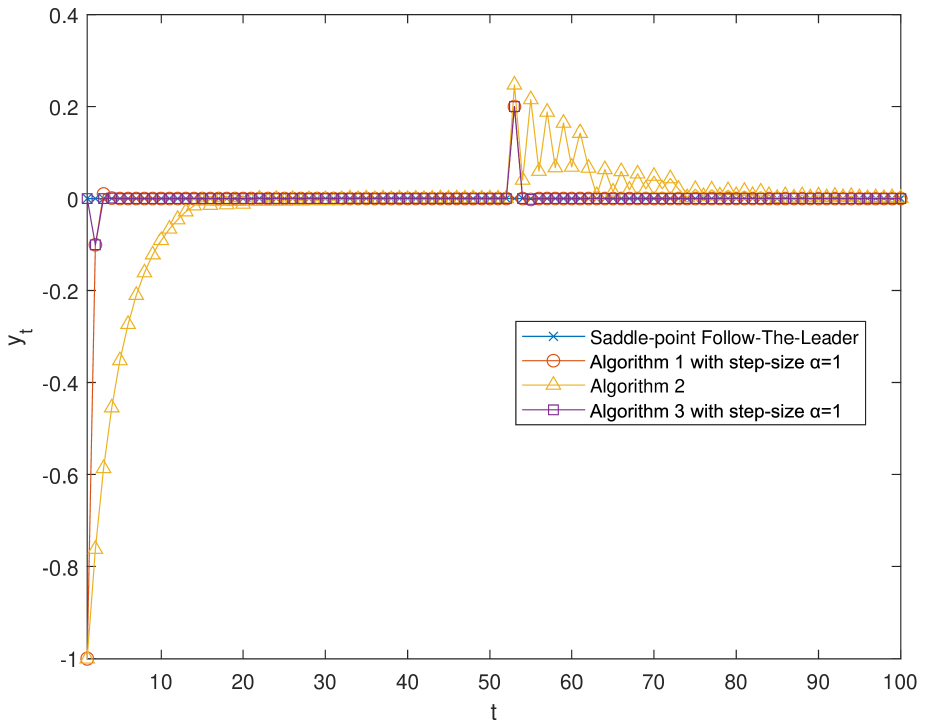}
\caption{Estimates of saddle-point Follow-The-Leader and Algorithms 1--3  with respect to $y$.}
\label{fig:y_first_order}
\end{figure}

\begin{figure}[htbp]
\centering
\includegraphics[width=0.95\linewidth]{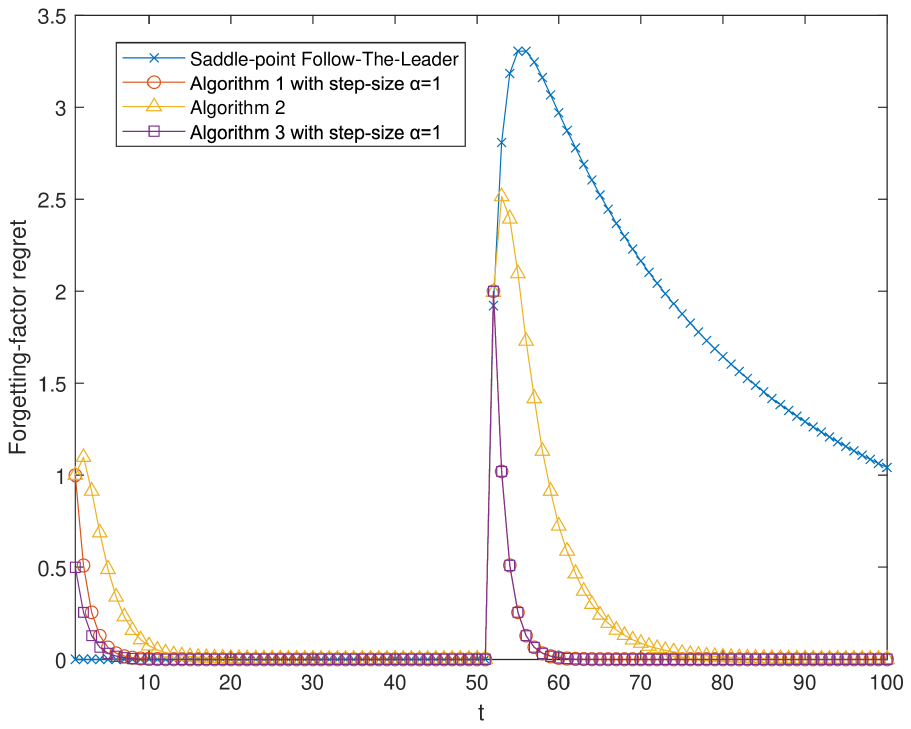}
\caption{Forgetting-factor regrets of saddle-point Follow-The-Leader and Algorithms 1--3.}
\label{fig:reg_first_order}
\end{figure}

Figs.~\ref{fig:x_first_order}--\ref{fig:y_first_order} show the estimates of
agents generated by Algorithms \ref{alg:ogda}--\ref{alg:zero} and Follow-The-Leader algorithm, respectively. 
Fig.~\ref{fig:reg_first_order} shows the forgetting-factor regrets of these  algorithms. The Follow-The-Leader algorithm has favorable historical-regret behavior, but it reacts slowly after the equilibrium switches from $(1,0)$ to $(-1,0)$. In contrast, the proposed online algorithms rapidly move toward the new NE and maintain small forgetting-factor regret. These results show that, for online zero-sum games, the forgetting-factor regret is an effective metric for evaluating the real-time tracking performance of NE-seeking algorithms.


\section{Conclusion}
\label{sec:conclusion}

This paper introduces a forgetting-factor regret for online zero-sum games and establishes its connection to terminal Nash equilibrium tracking. Three algorithms are  analyzed under this metric: projected online gradient descent-ascent, projection-free Frank-Wolfe, and deterministic-difference zeroth-order descent-ascent. The derived bounds show how tracking depends on equilibrium variation, payoff variation, cross-gradient coupling, and gradient-estimation accuracy. Sufficient conditions guaranteeing the proposed regret converges to zero are given. 
Future work includes weakening strong convexity-strong concavity, analyzing other online saddle-point algorithms under forgetting-factor regret, and developing adaptive choices of the forgetting factor and step sizes.

\appendices
\renewcommand{\thedefinition}{A\arabic{definition}}
\renewcommand{\thelemma}{A\arabic{lemma}}
\setcounter{definition}{0}
\setcounter{lemma}{0}
\section{Auxiliary Definitions and Lemmas}\label{Appendix:A}
\begin{definition}(Subgradient, strongly-convex, and convex-concave)\label{Def_A1}
\begin{itemize}
\item[a)] For a function $f(\cdot)$ of one variable, denote its domain by $\X \subset \mathbb{R}^d$. A vector-valued function $\partial f({x})\in \mathbb{R}^d$ is the subgradient of a non-smooth convex function $f(x):\X\to \mathbb{R}$ if for any ${x},{x}^{'}\in \X$,
$$f({x}^{'})-f({x})\geq \partial f({x})^\top ({x}^{'}-{x}).$$ $f(\cdot)$ is
said to be
$\sigma$-strongly convex on $\X$, if for any ${x},{x}^{'}\in \X$, $$f({x}^{'})-f({x})\geq \partial f({x})^\top ({x}^{'}-{x})+\frac{\sigma}{2}\|{x}^{'}-{x}\|^2,$$
and is said to be $L$-smooth on $\X$, if for any ${x},{x}^{'}\in \X$,
$$f({x}^{'})-f({x})\leq \partial f({x})^\top ({x}^{'}-{x})+\frac{L}{2}\|{x}^{'}-{x}\|^2.$$

\item[b)] For a function $f(\cdot,\cdot)$ of two variables, denote by $\X\times \Y\subset \mathbb{R}^{d_1}\times \mathbb{R}^{d_2}$ the domain of it. $f(\cdot,\cdot)$ is said to be convex-concave on $\X\times \Y$ if for any fixed ${y}$, $f({x},{y})$ is convex with respect to ${x}\in \X$ and for any fixed ${x}$, $f({x},{y})$ is concave with respect to ${y}\in \Y$. 
\end{itemize}
\end{definition}

\begin{lemma}\label{Lem:sigma}
Let $\mathcal{X}\subseteq\mathbb{R}^d$ be a convex set, and let
$f:\mathcal{X}\to\mathbb{R}$ be differentiable.

\begin{enumerate}
\item The following statements are equivalent:
\begin{enumerate}
    \item $f$ is $\sigma$-strongly convex on $\mathcal{X}$.
    \item For all ${x},{x}'\in\mathcal{X}$ and all
    $\lambda\in[0,1]$,
    \begin{align}
    f(\lambda {x}+(1-\lambda){x}')
    \le&
    \lambda f({x})+(1-\lambda)f({x}') \nonumber\\
    &-\frac{\sigma}{2}\lambda(1-\lambda)
    \|{x}-{x}'\|^2 .
    \end{align}
    \item For all ${x},{x}'\in\mathcal{X}$,
    \begin{align}
    \left(\nabla f({x})-\nabla f({x}')\right)^\top
    ({x}-{x}')
    \ge
    \sigma\|{x}-{x}'\|^2 .
    \end{align}
    \item The function
    $ f({x})-\frac{\sigma}{2}\|{x}\|^2$ 
    is convex on $\mathcal{X}$.
\end{enumerate}

\item If, in addition, $f$ is convex, then the following statements are equivalent:
\begin{enumerate}
    \item $f$ is $L$-smooth on $\mathcal{X}$.
    \item $\nabla f$ is Lipschitz continuous on $\mathcal{X}$ with constant $L$,
    i.e.,
    \begin{align}
    \|\nabla f({x})-\nabla f({x}')\|
    \le L\|{x}-{x}'\|,
    \quad
    \forall {x},{x}'\in\mathcal{X}.
    \end{align}
    \item The function $\frac{L}{2}\|{x}\|^2-f({x})$ 
    is convex on $\mathcal{X}$.
\end{enumerate}
\end{enumerate}
\end{lemma}

\begin{lemma}\label{LemA0}(see \cite{Kakutani1941,2013Auto})
If sets $\X$ and $\Y$ are compact and convex and the payoff function $f({x},{y})$ is continuous and convex-concave on $\X \times \Y$, then there exists a Nash equilibrium for the  zero-sum game with the payoff function $f({x},{y}):~\X \times \Y\rightarrow \mathbb{R}$.
\end{lemma}

\begin{lemma}\label{TLemA1}(see \cite{2012projection})
Assume that $f({x})$ is $\sigma$-strongly convex and differentiable over $\X$. For any ${x}\in \X$,%
\begin{align}
 f({x})- f({x}^{*})\geq \frac{\sigma}{2}\left\|{x}-{x}^{*}\right\|^{2},\label{zz1}
\end{align}
with ${x}^{*}\in \mathop{\arg\min}_{{x}\in \X}f({x}).$

\end{lemma}

\begin{lemma}\label{lem:projection}(see \cite{2010Distributed})
For the projection operator $\mathcal{P}_{\X}(\cdot):\mathbb{R}^d \to \X$,
\begin{equation*}
\|\mathcal{P}_{\X}({x})-\mathcal{P}_{\X}({x}^{'})\|\leq \|{x}-{x}^{'}\|,~~\forall {x},{x}^{'} \in \mathbb{R}^d.
\end{equation*}
\end{lemma}

\section{Proofs of  Lemmas in Section \ref{sec:problem}}\label{Appendix:B}
{\it{Proof of Lemma \ref{lem:exist_unique}: } }
Existence follows from Lemma \ref{LemA0}. To prove uniqueness, suppose that $(x_0^*,y_0^*)$ and $(x_1^*,y_1^*)$ are two NEs of $\mu_x$-strongly convex and $\mu_y$-strongly concave payoff $f$. By the saddle inequalities, we have
\begin{align}
f(x_0^*,y_0^*)&\le f(x_1^*,y_0^*)\le f(x_1^*,y_1^*),\nonumber\\
f(x_1^*,y_1^*)&\le f(x_0^*,y_1^*)\le f(x_0^*,y_0^*). \label{eq:unique_chain}
\end{align}
Hence
\begin{align}
f(x_0^*,y_0^*)=f(x_1^*,y_0^*)=
f(x_1^*,y_1^*)=f(x_0^*,y_1^*). \label{eq:unique_equal}
\end{align}
If $x_0^*\ne x_1^*$, then for every $\lambda\in(0,1)$, by Lemma \ref{Lem:sigma},  strong convexity of $f(\cdot,y_0^*)$ gives
\begin{align}
&f(\lambda x_0^*+(1-\lambda)x_1^*,y_0^*) \nonumber\\
&\le \lambda f(x_0^*,y_0^*)+(1-\lambda)f(x_1^*,y_0^*) -\frac{\mu_x}{2}\lambda(1-\lambda)\|x_0^*-x_1^*\|^2 \nonumber\\
&< f(x_0^*,y_0^*), \label{eq:unique_strict_x}
\end{align}\hfill $\blacksquare$
which contradicts the fact that $(x_0^*,y_0^*)$ is a NE. Thus $x_0^*=x_1^*$. The same argument applied to the $\mu_y$-strong convexity of $-f(x_0^*,\cdot)$ gives $y_0^*=y_1^*$.

{\it{Proof of Lemma \ref{lem:refined}: }}
By Assumption \ref{ass:smooth}, we have $L\ge \mu$. 
When $L=\mu$, the gradient difference satisfies $\nabla\phi(u)-\nabla\phi(v)=\mu(u-v)$, and the conclusion holds trivially. 

Next we consider the case $L>\mu$. Define
\begin{align}
g(u)\triangleq \phi(u)-\frac{L+\mu}{4}\|u\|^2.
\end{align}
Then, for all $u,v\in \X$, we have
\begin{align}
\nabla g(u)-\nabla g(v)
=\nabla\phi(u)-\nabla\phi(v)-\frac{L+\mu}{2}(u-v). \label{eq:g_grad}
\end{align}
By $\mu$-strong convexity and Lemma \ref{Lem:sigma}, it follows that 
\begin{align}
&(\nabla g(u)-\nabla g(v))^\top(u-v) \nonumber\\
&=(\nabla\phi(u)-\nabla\phi(v))^\top(u-v)
-\frac{L+\mu}{2}\|u-v\|^2 \nonumber\\
&\ge -\frac{L-\mu}{2}\|u-v\|^2. \label{eq:g_lower}
\end{align}
By $L$-smoothness and Lemma \ref{Lem:sigma}, we can obtain 
\begin{align}
&(\nabla g(u)-\nabla g(v))^\top(u-v) \nonumber\\
&=(\nabla\phi(u)-\nabla\phi(v))^\top(u-v)
-\frac{L+\mu}{2}\|u-v\|^2 \nonumber\\
&\le \frac{L-\mu}{2}\|u-v\|^2, \label{eq:g_upper}
\end{align}
and 
\begin{align}
\|\nabla g(u)-\nabla g(v)\|
&\le \|\nabla\phi(u)-\nabla\phi(v)\|
+\frac{L+\mu}{2}\|u-v\| \nonumber\\
&\le \frac{3L+\mu}{2}\|u-v\|. \label{eq:g_lip_rough}
\end{align}
Combining \eqref{eq:g_lower} and \eqref{eq:g_upper}, for all $u,v$, it follows that 
\begin{align}
\left|(\nabla g(u)-\nabla g(v))^\top(u-v)\right|
\le \frac{L-\mu}{2}\|u-v\|^2. \label{eq:g_inner_abs}
\end{align}
Consequently,  
\begin{align}
&\left|g(u)-g(v)-\nabla g(v)^\top(u-v)\right| \nonumber\\
&\le \int_0^1
\left|(\nabla g(v+\eta(u-v))-\nabla g(v))^\top(u-v)\right|d\eta \nonumber\\
&\le \int_0^1 \eta \frac{L-\mu}{2}\|u-v\|^2d\eta
=\frac{L-\mu}{4}\|u-v\|^2. \label{eq:g_remainder}
\end{align}

First, suppose that $u$ and $v$ are interior points of $\X$. By the convexity of $\X$, there exists a constant $\delta>0$ such that
\begin{align}
B_\delta(\theta u+(1-\theta)v)\subseteq\mathcal{X},\qquad \forall \theta\in[0,1].\label{B_delta}
\end{align}
Let
\begin{align}
N=\left\lceil
\frac{3L+\mu}{2(L-\mu)\delta}\|u-v\|
\right\rceil,
\end{align}
and 
\begin{align}
u_i=\frac{i}{N}u+\frac{N-i}{N}v.
\end{align}
Then \eqref{eq:g_lip_rough} gives
\begin{align}
\nonumber&\frac{1}{L-\mu}\|\nabla g(u_i)-\nabla g(u_{i-1})\|
\\ &\le \frac{3L+\mu}{2(L-\mu)}\|u_{i+1}-u_{i}\|=\frac{3L+\mu}{2(L-\mu)}\frac{\|u-v\|}{N}\le \delta. \label{eq:partition_radius}
\end{align}
Define
\begin{align}
d_i\triangleq \frac{1}{2}(u_i-u_{i-1})
-\frac{1}{L-\mu}(\nabla g(u_i)-\nabla g(u_{i-1})).
\end{align}
By \eqref{B_delta} and \eqref{eq:partition_radius}, we can obtain 
\begin{align}
u_i-d_i
&=\frac{1}{2}(u_i+u_{i-1})
+\frac{1}{L-\mu}(\nabla g(u_i)-\nabla g(u_{i-1})) \nonumber\\
&\in\mathcal{X},\nonumber\\
u_{i-1}+d_i
&=\frac{1}{2}(u_i+u_{i-1})
-\frac{1}{L-\mu}(\nabla g(u_i)-\nabla g(u_{i-1})) \nonumber\\
&\in\mathcal{X}. \label{eq:shifted_points}
\end{align}
Applying \eqref{eq:g_remainder} to the four admissible point pairs gives
\begin{align}
&g(u_{i-1}+d_i)-g(u_i)
-\nabla g(u_i)^\top(u_{i-1}+d_i-u_i) \nonumber\\
&\le \frac{L-\mu}{4}\|u_{i-1}+d_i-u_i\|^2,\nonumber\\
&g(u_i-d_i)-g(u_{i-1})
-\nabla g(u_{i-1})^\top(u_i-d_i-u_{i-1}) \nonumber\\
&\le \frac{L-\mu}{4}\|u_i-d_i-u_{i-1}\|^2,\nonumber\\
&-\big(g(u_{i-1}+d_i)-g(u_{i-1})
-\nabla g(u_{i-1})^\top d_i\big) \nonumber\\
&\le \frac{L-\mu}{4}\|d_i\|^2,\nonumber\\
&-\big(g(u_i-d_i)-g(u_i)+\nabla g(u_i)^\top d_i\big)
\le \frac{L-\mu}{4}\|d_i\|^2. \label{eq:four_remainders}
\end{align}
Adding \eqref{eq:four_remainders} yields
\begin{align}
&(\nabla g(u_i)-\nabla g(u_{i-1}))^\top(u_i-u_{i-1}-2d_i) \nonumber\\
&\le \frac{L-\mu}{2}\left(
\|u_{i-1}+d_i-u_i\|^2+\|d_i\|^2\right). \label{eq:sum_four}
\end{align}
Substituting the definition of $d_i$ into \eqref{eq:sum_four},
\begin{align}
&\frac{2}{L-\mu}\|\nabla g(u_i)-\nabla g(u_{i-1})\|^2 \nonumber\\
&\le \frac{L-\mu}{2}
\left\|\frac{1}{2}(u_i-u_{i-1}) +\frac{1}{L-\mu}(\nabla g(u_i)-\nabla g(u_{i-1}))\right\|^2 \nonumber\\&~+\frac{L-\mu}{2}
\left\|\frac{1}{2}(u_i-u_{i-1}) 
-\frac{1}{L-\mu}(\nabla g(u_i)-\nabla g(u_{i-1}))\right\|^2 \nonumber\\
&=\frac{L-\mu}{4}\|u_i-u_{i-1}\|^2 +\frac{1}{L-\mu}
\|\nabla g(u_i)-\nabla g(u_{i-1})\|^2.
\end{align}
Therefore,
\begin{align}
\|\nabla g(u_i)-\nabla g(u_{i-1})\|
\le \frac{L-\mu}{2}\|u_i-u_{i-1}\|.
\end{align}
Summing over the partition gives
\begin{align}
\|\nabla g(u)-\nabla g(v)\|
&\le \sum_{i=1}^{N}\|\nabla g(u_i)-\nabla g(u_{i-1})\| \nonumber\\
&\le \frac{L-\mu}{2}\sum_{i=1}^{N}\|u_i-u_{i-1}\| \nonumber\\
&=\frac{L-\mu}{2}\|u-v\|. \label{eq:g_int_bound}
\end{align}

For boundary points $u,v\in\mathcal{X}$, take interior points $\tilde u,\tilde v$ such that
$\|u-\tilde u\|\le\varepsilon$ and $\|v-\tilde v\|\le\varepsilon$. By \eqref{eq:g_lip_rough} and \eqref{eq:g_int_bound},
\begin{align}
&\|\nabla g(u)-\nabla g(v)\| \nonumber\\
&\le \|\nabla g(u)-\nabla g(\tilde u)\|
 \nonumber\\
&\quad+\|\nabla g(\tilde u)-\nabla g(\tilde v)\|
+\|\nabla g(v)-\nabla g(\tilde v)\| \nonumber\\
&\le \frac{L-\mu}{2}\big(\|u-v\|+2\varepsilon\big) \nonumber\\
&\quad +(3L+\mu)\varepsilon.
\end{align}
Letting $\varepsilon\downarrow0$ yields
\begin{align}
\|\nabla g(u)-\nabla g(v)\|\le \frac{L-\mu}{2}\|u-v\|.
\end{align}
Using \eqref{eq:g_grad} proves \eqref{eq:refined}.\hfill $\blacksquare$

\section{Proof of Theorem~\ref{thm:ogda}}\label{Appendix:C}

Fix $t$ and write
\begin{align}
a_t\triangleq x_t-x_t^*,\qquad b_t\triangleq y_t-y_t^*.
\end{align}
For compactness in the proof, define the line segments
\begin{align}
z_t^x(\eta)\triangleq x_t^*+\eta a_t,\qquad
z_t^y(\eta)\triangleq y_t^*+\eta b_t,\qquad \eta\in[0,1].
\end{align}
Define the cross-gradient perturbations
\begin{align}
d_t^x
&\triangleq \int_0^1
\big(\nabla_x f_t(z_t^x(\eta),y_t)
-\nabla_x f_t(z_t^x(\eta),y_t^*)\big)d\eta,\nonumber\\
d_t^y
&\triangleq -\int_0^1
\big(\nabla_y f_t(x_t,z_t^y(\eta))
-\nabla_y f_t(x_t^*,z_t^y(\eta))\big)d\eta . \label{eq:ogda_d_def}
\end{align}
Assumption~\ref{ass:cross} gives
\begin{align}
\|d_t^x\|\le L_\times\|b_t\|,\qquad
\|d_t^y\|\le L_\times\|a_t\|. \label{eq:d_bounds}
\end{align}
Moreover, by the fundamental theorem of calculus,
\begin{align}
&\langle d_t^x,a_t\rangle+\langle d_t^y,b_t\rangle \nonumber\\
&= \big(f_t(x_t,y_t)-f_t(x_t^*,y_t)\big)
-\big(f_t(x_t,y_t^*)-f_t(x_t^*,y_t^*)\big) \nonumber\\
&\quad-\big(f_t(x_t,y_t)-f_t(x_t,y_t^*)\big)
+\big(f_t(x_t^*,y_t)-f_t(x_t^*,y_t^*)\big) \nonumber\\
&=0, \label{eq:d_cancel}
\end{align}
where the first equality expands the two line integrals in \eqref{eq:ogda_d_def} along $x$ and $y$.

We next compare the true gradient at $(x_t,y_t)$ with the gradient that would be obtained if the opponent were at equilibrium. 
The interior-NE condition in Theorem~\ref{thm:ogda} gives
\begin{align}
\nabla_x f_t(x_t^*,y_t^*)=0,\qquad
\nabla_y f_t(x_t^*,y_t^*)=0. \label{eq:stationarity_used_ogda}
\end{align}
For the $x$-component, define
\begin{align}
A_t^x(\eta)&\triangleq \nabla_x f_t(x_t,y_t)
-\nabla_x f_t(z_t^x(\eta),y_t)
-\frac{L+\mu}{2}(1-\eta) a_t,\nonumber\\
B_t^x(\eta)&\triangleq \nabla_x f_t(z_t^x(\eta),y_t^*)
-\nabla_x f_t(x_t^*,y_t^*)
-\frac{L+\mu}{2}\eta a_t.
\end{align}
Using the definition of $d_t^x$ and the first identity in
\eqref{eq:stationarity_used_ogda}, we have the exact decomposition
\begin{align}
r_t^x-d_t^x-\frac{L+\mu}{2}a_t
=\int_0^1\big(A_t^x(\eta)+B_t^x(\eta)\big)d\eta . \label{eq:ogda_x_decomp}
\end{align}
Therefore, by the triangle inequality,
\begin{align}
\left\|r_t^x-d_t^x-\frac{L+\mu}{2}a_t\right\|
\le \int_0^1\|A_t^x(\eta)\|\,d\eta
+\int_0^1\|B_t^x(\eta)\|\,d\eta . \label{eq:ogda_x_triangle}
\end{align}
Lemma~\ref{lem:refined}, applied to $f_t(\cdot,y_t)$ and $f_t(\cdot,y_t^*)$, gives
\begin{align*}
\|A_t^x(\eta)\|&\le (1-\eta)\frac{L-\mu}{2}\|a_t\|,\\
\|B_t^x(\eta)\|&\le \eta\frac{L-\mu}{2}\|a_t\|.
\end{align*}
Substituting these two bounds into \eqref{eq:ogda_x_triangle} gives
\begin{align}
&\left\|r_t^x-d_t^x-\frac{L+\mu}{2}a_t\right\| \nonumber\\
&\le \frac{L-\mu}{2}\|a_t\|
\int_0^1(1-\eta)\,d\eta +\frac{L-\mu}{2}\|a_t\|
\int_0^1\eta\,d\eta \nonumber\\
&=\frac{L-\mu}{2}\|a_t\|. \label{eq:rd_x_full}
\end{align}

For the $y$-component, apply Lemma~\ref{lem:refined} to the strongly convex functions
$-f_t(x_t,\cdot)$ and $-f_t(x_t^*,\cdot)$. Define
\begin{align}
A_t^y(\eta)&\triangleq -\nabla_y f_t(x_t,y_t)
+\nabla_y f_t(x_t,z_t^y(\eta))
-\frac{L+\mu}{2}(1-\eta) b_t,\nonumber\\
B_t^y(\eta)&\triangleq -\nabla_y f_t(x_t^*,z_t^y(\eta))
+\nabla_y f_t(x_t^*,y_t^*)
-\frac{L+\mu}{2}\eta b_t .
\end{align}
Using the definition of $d_t^y$ and the second identity in
\eqref{eq:stationarity_used_ogda}, we similarly obtain
\begin{align}
-r_t^y-d_t^y-\frac{L+\mu}{2}b_t
=\int_0^1\big(A_t^y(\eta)+B_t^y(\eta)\big)d\eta . \label{eq:ogda_y_decomp}
\end{align}
The equality above is precisely the step where
$\nabla_y f_t(x_t^*,y_t^*)=0$ is needed. Lemma~\ref{lem:refined} gives
\begin{align*}
\|A_t^y(\eta)\|&\le (1-\eta)\frac{L-\mu}{2}\|b_t\|,\\
\|B_t^y(\eta)\|&\le \eta\frac{L-\mu}{2}\|b_t\|.
\end{align*}
Thus,
\begin{align}
\left\|-r_t^y-d_t^y-\frac{L+\mu}{2}b_t\right\|
\le \frac{L-\mu}{2}\|b_t\|. \label{eq:rd_y_full}
\end{align}
Let
\begin{align}
\kappa\triangleq 1-\alpha\frac{L+\mu}{2},\qquad
\tau\triangleq \alpha\frac{L-\mu}{2}.
\end{align}
By Algorithm \ref{alg:ogda}, Lemma \ref{lem:projection} and \eqref{eq:rd_x_full}, we have 
\begin{align*}
\|x_{t+1}-x_t^*\|
&\le \|x_t-\alpha r_t^x-x_t^*\|\\
&\le \left\|x_{t}-x_t^*-\alpha d_t^x-\alpha\frac{L+\mu}{2}(x_{t}-x_t^*))\right\|\\
&\quad+\alpha\left\|r_t^x-d_t^x-\frac{L+\mu}{2}(x_{t}-x_t^*))\right\|, 
\end{align*}
where the second inequality adds and subtracts $d_t^x+(L+\mu)(x_{t}-x_t^*)/2$. Hence, by \eqref{eq:rd_x_full},
\begin{align}
\|x_{t+1}-x_t^*\|
\le \|\kappa (x_{t}-x_t^*)-\alpha d_t^x\|+\tau\|x_{t}-x_t^*\|. \label{eq:dxx_full}
\end{align}
Likewise,
\begin{align}
\|y_{t+1}-y_t^*\|
\le \|\kappa (y_{t}-y_t^*)-\alpha d_t^y\|+\tau\|y_{t}-y_t^*\|. \label{eq:dyy_full}
\end{align}
From \eqref{eq:d_bounds} and \eqref{eq:d_cancel},
\begin{align}
\nonumber&\|\kappa (x_{t}-x_t^*)-\alpha d_t^x\|^2
+\|\kappa (y_{t}-y_t^*)-\alpha d_t^y\|^2\\
\nonumber&=\kappa^2(\|x_{t}-x_t^*\|^2+\|y_{t}-y_t^*\|^2)
-2\alpha\kappa\big(\langle d_t^x,x_{t}-x_t^*\rangle\\
\nonumber&\quad+\langle d_t^y,y_{t}-y_t^*\rangle\big)+\alpha^2(\|d_t^x\|^2+\|d_t^y\|^2)
\\
&\le \left(\kappa^2+\alpha^2L_\times^2\right)
\left(\|x_{t}-x_t^*\|^2+\|y_{t}-y_t^*\|^2\right). \label{eq:sq_full}
\end{align}
By Cauchy's inequality and \eqref{eq:sq_full},
\begin{align}
&\|\kappa (x_{t}\!-\!x_t^*)\!-\!\alpha d_t^x\|\|x_{t}\!-\!x_t^*\|
+\|\kappa (y_{t}\!-\!y_t^*)\!-\!\alpha d_t^y\|\|y_{t}\!-\!y_t^*\| \nonumber\\
&\le \sqrt{\kappa^2+\alpha^2L_\times^2}
\left(\|x_{t}-x_t^*\|^2+\|y_{t}-y_t^*\|^2\right). \label{eq:cr_full}
\end{align}
Squaring \eqref{eq:dxx_full} and \eqref{eq:dyy_full} and summing gives
\begin{align*}
&\|x_{t+1}-x_t^*\|^2+\|y_{t+1}-y_t^*\|^2\\\nonumber
&\le
\left(\|\kappa (x_{t}-x_t^*)-\alpha d_t^x\|+\tau\|x_{t}-x_t^*\|\right)^2
\\
&\quad+\left(\|\kappa (y_{t}-y_t^*)-\alpha d_t^y\|+\tau\|y_{t}-y_t^*\|\right)^2 .
\end{align*}
Combining this estimate with \eqref{eq:sq_full}--\eqref{eq:cr_full} yields
\begin{align}
&\|x_{t+1}-x_t^*\|^2+\|y_{t+1}-y_t^*\|^2 \nonumber\\
\le& \!\Big(\kappa^2\!\!+\!\alpha^2L_\times^2
\!\!+\!\!2\tau\sqrt{\kappa^2\!+\!\alpha^2L_\times^2}\!+\!\tau^2\Big)
\left(\|x_{t}\!-\!x_t^*\|^2\!\!+\!\!\|y_{t}\!-\!y_t^*\|^2\right) \nonumber\\
=&\beta\left(\|x_{t}-x_t^*\|^2+\|y_{t}-y_t^*\|^2\right). \label{eq:fixed_ne_contraction}
\end{align}
The last equality is the definition of $\beta$ in \eqref{eq:beta_gda}. The range
\eqref{eq:alpha_range} is equivalent to $\beta<1$, and $\beta\ge0$ follows from
\eqref{eq:fixed_ne_contraction}.

For $t\ge1$, the temporal variation of the equilibrium and compactness give
\begin{align*}
&\|x_t-x_t^*\|^2+\|y_t-y_t^*\|^2\\
&\le \|x_t-x_{t-1}^*\|^2+2D_\X V_t^x+(V_t^x)^2\\
&\quad+\|y_t-y_{t-1}^*\|^2+2D_\Y V_t^y+(V_t^y)^2\\
&\le \beta\big(\|x_{t-1}-x_{t-1}^*\|^2
+\|y_{t-1}-y_{t-1}^*\|^2\big)\\
&\quad+2D_\X V_t^x+2D_\Y V_t^y+(V_t^x)^2+(V_t^y)^2 .
\end{align*}
Here the compactness bounds use the diameters in \eqref{eq:diameters}, and the
last inequality applies \eqref{eq:fixed_ne_contraction} at time $t-1$.
Iterating the above recursion gives
\begin{align}
&\|x_t-x_t^*\|^2+\|y_t-y_t^*\|^2 \nonumber\\
&\le \beta^{t}\big(\|x_0-x_0^*\|^2+\|y_0-y_0^*\|^2\big) \nonumber\\
&\quad+
\sum_{i=1}^{t}\beta^{t-i}
\big(2D_\X V_i^x+2D_\Y V_i^y \nonumber\\
&\quad+(V_i^x)^2+(V_i^y)^2\big), 
\label{eq:ogda_compact_iter_full}
\end{align}
which yields \eqref{eq:ogda_compact_e}. 
By Assumption \ref{ass:smooth}, Lemma \ref{Lem:sigma}  and the stationarity conditions in \eqref{eq:stationarity_used_ogda}, we have 
\begin{align}
&f_t(x_t,y_t^*)-f_t(x_t^*,y_t) \nonumber\\
&=f_t(x_t,y_t^*)-f_t(x_t^*,y_t^*) \nonumber\\
&\quad+f_t(x_t^*,y_t^*)-f_t(x_t^*,y_t) \nonumber\\
&\le \frac{L_x}{2}\|x_t-x_t^*\|^2
+\frac{L_y}{2}\|y_t-y_t^*\|^2 \nonumber\\
&\le \frac{L}{2}\left(\|x_t-x_t^*\|^2+\|y_t-y_t^*\|^2\right).
\label{eq:smooth_gap_bound}
\end{align}
Combining \eqref{eq:smooth_gap_bound} with \eqref{eq:ogda_compact_iter_full},
and then multiplying by $\rho^{T-t}$ and summing over $t$, gives
\begin{align}
\RegF
&\le \frac{L}{2}\sum_{t=1}^{T}\rho^{T-t}\beta^{t}
\big(\|x_0-x_0^*\|^2+\|y_0-y_0^*\|^2\big) \nonumber\\
&\quad+\frac{L}{2}\sum_{i=1}^{T}\sum_{t=i}^{T}
\rho^{T-t}\beta^{t-i}
\big(2D_\X V_i^x+2D_\Y V_i^y \nonumber\\
&\quad+(V_i^x)^2+(V_i^y)^2\big),
\end{align}
which is \eqref{eq:ogda_compact_regret}. The convergence statement follows by
applying the standard fact that a vanishing sequence convolved with a stable
geometric kernel remains vanishing.\hfill $\blacksquare$

\section{Proof of Theorem~\ref{thm:fw}}\label{Appendix:D}

By the Frank-Wolfe update in Algorithm \ref{alg:fw}, the smoothness of $f_t(\cdot,y_t^*)$, Assumption \ref{ass:set} and \ref{ass:cross}, we can obtain 
\begin{align}
&f_t(x_{t+1},y_t^*)-f_t(x_t^*,y_t^*) \nonumber\\
&=f_t((1-\eta_t)x_t+\eta_t s_t^x,y_t^*)
-f_t(x_t^*,y_t^*) \nonumber\\
&\le f_t(x_t,y_t^*)-f_t(x_t^*,y_t^*) \nonumber\\
&+\eta_t\big\langle\nabla_x f_t(x_t,y_t^*), s_t^x-x_t\big\rangle \nonumber\\
&\quad+\eta_t^2\frac{L_x}{2}\|s_t^x-x_t\|^2 \nonumber\\
&\le f_t(x_t,y_t^*)-f_t(x_t^*,y_t^*) +\eta_t\langle r_t^x,s_t^x-x_t\rangle \nonumber\\
&\quad+\eta_tD_\X L_{xy}\|y_t-y_t^*\| +\eta_t^2\frac{L_xD_\X^2}{2} \nonumber\\
&\le f_t(x_t,y_t^*)-f_t(x_t^*,y_t^*)-\eta_t g_t^x
 \nonumber\\
&\quad+\eta_t D_\X L_{xy}\|y_t-y_t^*\| +\eta_t^2\frac{L_xD_\X^2}{2}. \label{eq:fw_x_step}
\end{align}
Similarly, we can get 
\begin{align}
&f_t(x_t^*,y_t^*)-f_t(x_t^*,y_{t+1}) \nonumber\\
&\le f_t(x_t^*,y_t^*)-f_t(x_t^*,y_t)-\eta_t g_t^y \nonumber\\
&\quad+\eta_t D_\Y L_{yx}\|x_t-x_t^*\| +\eta_t^2\frac{L_yD_\Y^2}{2}. \label{eq:fw_y_step}
\end{align}
Adding \eqref{eq:fw_x_step} and \eqref{eq:fw_y_step}, combining \eqref{def:C}, Lemma \ref{TLemA1} yields
\begin{align}
&f_t(x_{t+1},y_t^*)-f_t(x_t^*,y_{t+1}) \nonumber\\
&\le f_t(x_t,y_t^*)-f_t(x_t^*,y_t)-\eta_t g_t +\eta_t D_\X L_{xy}\|y_t-y_t^*\| \nonumber\\
&\quad+\eta_t D_\Y L_{yx}\|x_t-x_t^*\| +\eta_t^2\frac{C}{2} \nonumber\\
&\le f_t(x_t,y_t^*)-f_t(x_t^*,y_t)-\eta_t g_t \nonumber\\
&\quad+\eta_t D_\X L_{xy}\sqrt{{2}/{\mu_y}} \sqrt{f_t(x_t^*,y_t^*)-f_t(x_t^*,y_t)} \nonumber\\
&\quad+\eta_t D_\Y L_{yx}\sqrt{{2}/{\mu_x}} \sqrt{f_t(x_t,y_t^*)-f_t(x_t^*,y_t^*)} +\eta_t^2\frac{C}{2} \nonumber\\
&\le f_t(x_t,y_t^*)-f_t(x_t^*,y_t)-\eta_t g_t \nonumber\\
&\quad+\eta_t M_c
\sqrt{2\big(f_t(x_t,y_t^*)-f_t(x_t^*,y_t^*)\big)}
\nonumber\\
&\quad+\eta_t M_c
\sqrt{2\big(f_t(x_t^*,y_t^*)-f_t(x_t^*,y_t)\big)}
+\eta_t^2\frac{C}{2} \nonumber\\
&\le f_t(x_t,y_t^*)-f_t(x_t^*,y_t)-\eta_t g_t \nonumber\\
&\quad+2\eta_t M_c
\sqrt{f_t(x_t,y_t^*)-f_t(x_t^*,y_t)} +\eta_t^2\frac{C}{2},  \label{eq:fw_basic}
\end{align}
where the last two
inequalities 
$\sqrt{a}+\sqrt{b}\le\sqrt{2(a+b)}$ for $a,b\ge0$ and the notations
\begin{align}
M_c&\triangleq \max\left\{\frac{D_\X L_{xy}}{\sqrt{\mu_y}},
\frac{D_\Y L_{yx}}{\sqrt{\mu_x}}\right\},\\
S_c&\triangleq \min\left\{\sqrt{\mu_y}\delta_\Y,\sqrt{\mu_x}\delta_\X\right\}.\label{def:SC}
\end{align}

Since $x_t^*$ is at least $\delta_\X$ away from $\partial\X$ by \eqref{def:delta}, we have  
$x_t^*-\delta_\X r_t^x/\|r_t^x\|\in\X$. The definition of $s_t^x$ gives
\begin{align}
g_t^x
&=\langle r_t^x,x_t-s_t^x\rangle \nonumber\\
&\ge \left\langle r_t^x,
x_t-x_t^*+\delta_\X\frac{r_t^x}{\|r_t^x\|}\right\rangle \nonumber\\
&=\delta_\X\|r_t^x\|+\langle r_t^x,x_t-x_t^*\rangle . \label{eq:fw_rx}
\end{align}
Similarly, since $y_t^*+\delta_\Y r_t^y/\|r_t^y\|\in\Y$ and $s_t^y$
maximizes the linear oracle over $\Y$, we can get 
\begin{align}
g_t^y
\ge \delta_\Y\|r_t^y\|+\langle r_t^y,y_t^*-y_t\rangle . \label{eq:fw_ry}
\end{align}
By strong convexity in $x$ and strong concavity in $y$,
\begin{align}
&\langle r_t^x,x_t-x_t^*\rangle
+\langle r_t^y,y_t^*-y_t\rangle \nonumber\\
&\ge f_t(x_t,y_t)-f_t(x_t^*,y_t)
+\frac{\mu_x}{2}\|x_t-x_t^*\|^2 \nonumber\\
&\quad+f_t(x_t,y_t^*)-f_t(x_t,y_t)
+\frac{\mu_y}{2}\|y_t-y_t^*\|^2 \nonumber\\
&=f_t(x_t,y_t^*)-f_t(x_t^*,y_t)
+\frac{\mu_x}{2}\|x_t-x_t^*\|^2 \nonumber\\
&\quad+\frac{\mu_y}{2}\|y_t-y_t^*\|^2
\ge0. \label{eq:fw_monotone}
\end{align}
Adding \eqref{eq:fw_rx} and \eqref{eq:fw_ry}, and then using
\eqref{eq:fw_monotone}, yields
\begin{align}
g_t\ge \delta_\X\|r_t^x\|+\delta_\Y\|r_t^y\|. \label{eq:fw_g_lower}
\end{align}
Using strong concavity in $y$ and strong convexity in $x$, we obtain
\begin{align*}
&f_t(x_t,y_t^*)-f_t(x_t^*,y_t)\\
&\le \langle r_t^y,y_t^*-y_t\rangle
-\frac{\mu_y}{2}\|y_t^*-y_t\|^2\\
&\quad+\langle r_t^x,x_t-x_t^*\rangle
-\frac{\mu_x}{2}\|x_t-x_t^*\|^2.
\end{align*}
By Cauchy's inequality and $\max_{u\ge0}\{au-\mu u^2/2\}=a^2/(2\mu)$,
\begin{align*}
f_t(x_t,y_t^*)-f_t(x_t^*,y_t)
\le \frac{\|r_t^y\|^2}{2\mu_y}
+\frac{\|r_t^x\|^2}{2\mu_x}.
\end{align*}
Using \eqref{def:SC} and \eqref{eq:fw_g_lower} gives
\begin{align}
&f_t(x_t,y_t^*)-f_t(x_t^*,y_t) \nonumber\\
&\le \frac{\delta_\Y^2\|r_t^y\|^2+\delta_\X^2\|r_t^x\|^2}{2S_c^2}
\le \frac{g_t^2}{2S_c^2}. \label{eq:fw_gap_gt}
\end{align}
Substituting \eqref{eq:fw_gap_gt} into \eqref{eq:fw_basic} and using the definition of $\nu$ gives
\begin{align}
&f_t(x_{t+1},y_t^*)-f_t(x_t^*,y_{t+1}) \nonumber\\
&\le f_t(x_t,y_t^*)-f_t(x_t^*,y_t)-\eta_t g_t
+\sqrt{2}\eta_t\frac{M_c}{S_c}g_t
+\eta_t^2\frac{C}{2} \nonumber\\
&= f_t(x_t,y_t^*)-f_t(x_t^*,y_t)-\eta_t\nu g_t
+\eta_t^2\frac{C}{2}. \label{eq:fw_descent}
\end{align}
Furthermore, the interior stationarity conditions
$\nabla_x f_t(x_t^*,y_t^*)=0$ and $\nabla_y f_t(x_t^*,y_t^*)=0$,
Assumptions~\ref{ass:smooth} and~\ref{ass:cross}, and the diameter bounds in
\eqref{eq:diameters} imply
\begin{align}
g_t
&=g_t^x+g_t^y \nonumber\\
&\le \|r_t^x\|\,\|x_t-s_t^x\|
+\|r_t^y\|\,\|s_t^y-y_t\| \nonumber\\
&\le \big(L_x\|x_t-x_t^*\|+L_{xy}\|y_t-y_t^*\|\big)D_\X \nonumber\\
&\quad+\big(L_y\|y_t-y_t^*\|+L_{yx}\|x_t-x_t^*\|\big)D_\Y \nonumber\\
&\le L_xD_\X^2+L_yD_\Y^2+(L_{xy}+L_{yx})D_\X D_\Y=C. \label{eq:fw_stepsize_feasible}
\end{align}
Since $\nu>0$ implies $\nu<1$, \eqref{eq:fw_stepsize_feasible} gives
$\nu g_t/C\le1$.
Therefore $\eta_t=\nu g_t/C\in[0,1]$ and minimizes the right-hand side of
\eqref{eq:fw_descent}. Combining \eqref{eq:fw_descent} and \eqref{eq:fw_gap_gt},
\begin{align}
&f_t(x_{t+1},y_t^*)-f_t(x_t^*,y_{t+1}) \nonumber\\
&\le f_t(x_t,y_t^*)-f_t(x_t^*,y_t)-\frac{\nu^2g_t^2}{2C} \nonumber\\
&\le f_t(x_t,y_t^*)-f_t(x_t^*,y_t) \nonumber\\
&\quad-\frac{\nu^2\min\{\mu_y\delta_\Y^2,\mu_x\delta_\X^2\}}{C}
\big(f_t(x_t,y_t^*)-f_t(x_t^*,y_t)\big) \nonumber\\
&=\beta_{\rm FW}\big(f_t(x_t,y_t^*)-f_t(x_t^*,y_t)\big). \label{eq:fw_one_step}
\end{align}
Here $S_c^2=\min\{\mu_y\delta_\Y^2,\mu_x\delta_\X^2\}$ and the last equality is the definition of $\beta_{\rm FW}$. Since $L_x\ge\mu_x$, $L_y\ge\mu_y$, $\delta_\X\le D_\X$, $\delta_\Y\le D_\Y$, and $\nu\in(0,1)$, we have $S_c^2\nu^2<C$. Hence $\beta_{\rm FW}\in(0,1)$.

We now compare consecutive stages. By definition of $F_t^{\sup}$,
\begin{align*}
&f_t(x_t,y_t^*)-f_t(x_t^*,y_t)\\
&\le f_{t-1}(x_t,y_t^*)-f_{t-1}(x_t^*,y_t)+2F_t^{\sup}.
\end{align*}
Adding and subtracting the previous NE terms gives
\begin{align*}
&f_t(x_t,y_t^*)-f_t(x_t^*,y_t)\\
&\le f_{t-1}(x_t,y_{t-1}^*)-f_{t-1}(x_{t-1}^*,y_t)\\
&\quad+\big(f_{t-1}(x_t,y_t^*)-f_{t-1}(x_t,y_{t-1}^*)\big)\\
&\quad+\big(f_{t-1}(x_{t-1}^*,y_t)-f_{t-1}(x_t^*,y_t)\big)+2F_t^{\sup}.
\end{align*}
We now bound the two variation terms explicitly. Since $(x_{t-1}^*,y_{t-1}^*)$
is an interior NE, stationarity gives
\begin{align}
\nabla_x f_{t-1}(x_{t-1}^*,y_{t-1}^*)=0,\qquad
\nabla_y f_{t-1}(x_{t-1}^*,y_{t-1}^*)=0. \label{eq:fw_prev_stationarity}
\end{align}
By $L_y$-smoothness in $y$, Assumption~\ref{ass:cross}, and
\eqref{eq:fw_prev_stationarity},
\begin{align}
\|\nabla_y f_{t-1}(x_t,y_{t-1}^*)\|
&\le L_{yx}\|x_t-x_{t-1}^*\|
\le D_\X L_{yx}, \nonumber\\
\|\nabla_x f_{t-1}(x_{t-1}^*,y_t)\|
&\le L_{xy}\|y_t-y_{t-1}^*\|
\le D_\Y L_{xy}. \label{eq:fw_variation_grad_bounds}
\end{align}
Therefore,
\begin{align}
&f_{t-1}(x_t,y_t^*)-f_{t-1}(x_t,y_{t-1}^*) \nonumber\\
&\le
\left\langle\nabla_y f_{t-1}(x_t,y_{t-1}^*),
y_t^*-y_{t-1}^*\right\rangle
+\frac{L_y}{2}(V_t^y)^2 \nonumber\\
&\le D_\X L_{yx}V_t^y+\frac{L_y}{2}(V_t^y)^2. \label{eq:fw_var_y}
\end{align}
Similarly, by $L_x$-smoothness in $x$ and \eqref{eq:fw_variation_grad_bounds},
\begin{align}
&f_{t-1}(x_{t-1}^*,y_t)-f_{t-1}(x_t^*,y_t) \nonumber\\
&\le
\left\langle\nabla_x f_{t-1}(x_{t-1}^*,y_t),
x_{t-1}^*-x_t^*\right\rangle
+\frac{L_x}{2}(V_t^x)^2 \nonumber\\
&\le D_\Y L_{xy}V_t^x+\frac{L_x}{2}(V_t^x)^2. \label{eq:fw_var_x}
\end{align}
Substituting \eqref{eq:fw_var_y}--\eqref{eq:fw_var_x} into the preceding
display yields
\begin{align*}
&f_t(x_t,y_t^*)-f_t(x_t^*,y_t)\\
&\le f_{t-1}(x_t,y_{t-1}^*)-f_{t-1}(x_{t-1}^*,y_t)\\
&\quad+D_\X L_{yx}V_t^y+D_\Y L_{xy}V_t^x
+\frac{L_y}{2}(V_t^y)^2\\
&\quad+\frac{L_x}{2}(V_t^x)^2+2F_t^{\sup}.
\end{align*}
Using the one-step contraction \eqref{eq:fw_one_step} at time $t-1$, we obtain
\begin{align*}
&f_t(x_t,y_t^*)-f_t(x_t^*,y_t)\\
&\le \beta_{\rm FW}
\big(f_{t-1}(x_{t-1},y_{t-1}^*)
-f_{t-1}(x_{t-1}^*,y_{t-1})\big)\\
&\quad+D_\X L_{yx}V_t^y+D_\Y L_{xy}V_t^x\\
&\quad+\frac{L_y}{2}(V_t^y)^2+\frac{L_x}{2}(V_t^x)^2
+2F_t^{\sup}.
\end{align*}
Iterating this scalar recursion yields
\begin{align}
&f_t(x_t,y_t^*)-f_t(x_t^*,y_t) \nonumber\\
&\le \beta_{\rm FW}^{t}
\big(f_0(x_0,y_0^*)-f_0(x_0^*,y_0)\big) \nonumber\\
&\quad+\sum_{i=1}^{t}\beta_{\rm FW}^{t-i}
\Big(D_\X L_{yx}V_i^y+D_\Y L_{xy}V_i^x \nonumber\\
&\quad+\frac{L_y}{2}(V_i^y)^2+\frac{L_x}{2}(V_i^x)^2
+2F_i^{\sup}\Big). \label{eq:fw_iter_full}
\end{align}
Finally,
\begin{align}
&\RegF
=\sum_{t=1}^{T}\rho^{T-t}
\left(f_t(x_t,y_t^*)-f_t(x_t^*,y_t)\right)
 \nonumber\\
&\le \sum_{t=1}^{T}\rho^{T-t}\beta_{\rm FW}^{t}
\big(f_0(x_0,y_0^*)-f_0(x_0^*,y_0)\big) \nonumber\\
&\quad+\sum_{t=1}^{T}\rho^{T-t}\sum_{i=1}^{t}
\beta_{\rm FW}^{t-i} \nonumber\\
&\quad\cdot\Big(D_\X L_{yx}V_i^y+D_\Y L_{xy}V_i^x
+\frac{L_y}{2}(V_i^y)^2 \nonumber\\
&\quad+\frac{L_x}{2}(V_i^x)^2+2F_i^{\sup}\Big) \nonumber\\
&=\sum_{t=1}^{T}\rho^{T-t}\beta_{\rm FW}^{t}
\big(f_0(x_0,y_0^*)-f_0(x_0^*,y_0)\big) \nonumber\\
&\quad+\sum_{i=1}^{T}\sum_{t=i}^{T}
\rho^{T-t}\beta_{\rm FW}^{t-i}
\Big(D_\X L_{yx}V_i^y+D_\Y L_{xy}V_i^x \nonumber\\
&\quad+\frac{L_y}{2}(V_i^y)^2+\frac{L_x}{2}(V_i^x)^2 +2F_i^{\sup}\Big),
\end{align}
where the last equality exchanges the order of summation.
This proves \eqref{eq:fw_gap_bound} and \eqref{eq:fw_regret}. The convergence
statement again follows from geometric convolution.\hfill $\blacksquare$

\section{Proof of Theorem~\ref{thm:zero}}\label{Appendix:E}

Let $r_t^x=\nabla_x f_t(x_t,y_t)$ and $r_t^y=\nabla_y f_t(x_t,y_t)$. By the Lagrange mean-value theorem, for each coordinate $k$ there exists $\theta_k\in(-\delta_t,\delta_t)$ such that
\begin{align}
h_t^x[k]
&=\frac{f_t(x_t+\delta_t e_x^k,y_t)-f_t(x_t-\delta_t e_x^k,y_t)}
{2\delta_t} \nonumber\\
&=\nabla_x f_t(x_t+\theta_k e_x^k,y_t)[k], 
\label{eq:zero_mvt_x}
\end{align}
Therefore, by Assumption~\ref{ass:smooth},
\begin{align}
\|h_t^x-r_t^x\|
&\le \sum_{k=1}^{n_x}
\left|\nabla_xf_t(x_t+\theta_ke_x^k,y_t)
-\nabla_xf_t(x_t,y_t)\right| \nonumber\\
&\le L_x\sum_{k=1}^{n_x}|\theta_k|
\le n_x\delta_t L_x
\le c_t. \label{eq:zero_grad_error_x}
\end{align}
The last inequality follows from the choice of $\delta_t$ in \eqref{eq:delta_choice}.
Similarly,
\begin{align}
\|h_t^y-r_t^y\|\le c_t. \label{eq:zero_grad_error_y}
\end{align}
The proof of \eqref{eq:zero_grad_error_y} is identical, using the $y$-coordinate central differences and the constant $L_y$.
Then, by the update rule in Algorithm \ref{alg:zero},  
\begin{align*}
&\|x_{t+1}-x_{t+1}^*\|^2\\
&=\|\mathcal{P}_{\X}\big(x_t-\alpha r_t^x-\alpha(h_t^x-r_t^x)\big)-x_t^*+x_t^*-x_{t+1}^*\|^2\\
&\le \|\mathcal{P}_{\X}\big(x_t\!-\!\alpha r_t^x\!-\!\alpha(h_t^x\!-\!r_t^x)\big)\!-\!x_t^*\|^2
\!+\!2D_\X V_{t+1}^x\!+\!(V_{t+1}^x)^2.
\end{align*}
Moreover, by Lemma \ref{lem:projection} and \eqref{eq:zero_grad_error_x},
\begin{align*}
&\|\mathcal{P}_{\X}\big(x_t-\alpha r_t^x-\alpha(h_t^x-r_t^x)\big)-x_t^*\|^2
\\\nonumber&=\|\mathcal{P}_{\X}(x_t-\alpha r_t^x)-x_t^*+\mathcal{P}_{\X}\big(x_t-\alpha r_t^x-\alpha(h_t^x-r_t^x)\big)
\\\nonumber&\quad  -\mathcal{P}_{\X}\big(x_t-\alpha r_t^x\big)\|^2\\
&\le \|\mathcal{P}_{\X}(x_t-\alpha r_t^x)-x_t^*\|^2
+2D_\X\alpha c_t+\alpha^2c_t^2\\
&\le \|x_t-\alpha r_t^x-x_t^*\|^2
+2D_\X\alpha c_t+\alpha^2c_t^2 .
\end{align*}
Combining the two estimates leads to
\begin{align}
\nonumber\|x_{t+1}-x_{t+1}^*\|^2
&\le \|x_t-\alpha r_t^x-x_t^*\|^2 +2D_\X(\alpha c_t+V_{t+1}^x)
\\
&\quad+\alpha^2c_t^2+(V_{t+1}^x)^2. \label{eq:zero_x_compact}
\end{align}
Similarly,
\begin{align}
\|y_{t+1}-y_{t+1}^*\|^2
&\le \|y_t+\alpha r_t^y-y_t^*\|^2
+2D_\Y(\alpha c_t+V_{t+1}^y) \nonumber\\
&\quad+\alpha^2c_t^2+(V_{t+1}^y)^2. \label{eq:zero_y_compact}
\end{align}
Combining \eqref{eq:zero_x_compact}--\eqref{eq:zero_y_compact} with
\eqref{eq:fixed_ne_contraction} gives
\begin{align}
&\|x_{t+1}-x_{t+1}^*\|^2+\|y_{t+1}-y_{t+1}^*\|^2 \nonumber\\
&\le \beta\big(\|x_t-x_t^*\|^2+\|y_t-y_t^*\|^2\big) \nonumber\\
&\quad+2D_\X(\alpha c_t+V_{t+1}^x) +2D_\Y(\alpha c_t+V_{t+1}^y) \nonumber\\
&\quad+2\alpha^2c_t^2+(V_{t+1}^x)^2 +(V_{t+1}^y)^2.
\end{align}
Here \eqref{eq:fixed_ne_contraction} is applied to the exact-gradient projected
descent-ascent step at time $t$. 
Thus, for $t\ge0$,
\begin{align}
&\|x_t-x_t^*\|^2+\|y_t-y_t^*\|^2 \nonumber\\
&\le \beta^{t}\big(\|x_0-x_0^*\|^2+\|y_0-y_0^*\|^2\big) \nonumber\\
&\quad+\sum_{i=1}^{t}\beta^{t-i}
\Big(2D_\X(\alpha c_{i-1}+V_i^x)
+2D_\Y(\alpha c_{i-1}+V_i^y) \nonumber\\
&\quad+2\alpha^2c_{i-1}^2+(V_i^x)^2 +(V_i^y)^2\Big), \label{eq:zero_compact_iter_full}
\end{align}
which follows by iterating the compact-set recursion. This is
\eqref{eq:zero_compact_e}. Finally,
\begin{align}
&\RegF
=\sum_{t=1}^{T}\rho^{T-t}
\left(f_t(x_t,y_t^*)-f_t(x_t^*,y_t)\right) \nonumber\\
&\le \frac{L}{2}\sum_{t=1}^{T}\rho^{T-t}\beta^{t}
\big(\|x_0-x_0^*\|^2+\|y_0-y_0^*\|^2\big) \nonumber\\
&\nonumber\quad+\frac{L}{2}\sum_{i=1}^{T}\sum_{t=i}^{T}
\rho^{T-t}\beta^{t-i} \left(2D_\X(\alpha c_{i-1}+V_i^x)\right. 
\\&\left.\quad\!+\!2D_\Y(\alpha c_{i-1}\!+\!V_i^y)\!+\!2\alpha^2c_{i-1}^2\!+\!(V_i^x)^2 \!+\!(V_i^y)^2\right),
\end{align}
where the last inequality uses \eqref{eq:smooth_gap_bound}, substitutes
\eqref{eq:zero_compact_iter_full}, and exchanges the summations. This is
\eqref{eq:zero_compact_regret}. The convergence statement follows from the same
geometric-kernel argument used in Theorem~\ref{thm:ogda}. \hfill $\blacksquare$

\bibliographystyle{ieeetr}
\bibliography{reference}

@book{2009Ben,
    title={Robust optimization},
    author={Ben-Tal, A. and El Ghaoui, L.  and Nemirovski, A. },
    year={2009},
    publisher={Princeton University Press},
}

@inproceedings{1996Freund,
  title={Game theory, online prediction and boosting},
  author={Freund, Y. and Schapire, R. E. },
  booktitle={Proceedings of the 9th Annual Conference on Computational Learning Theory},
  pages={325–332},
  year={1996},
}

@inproceedings{2014Good,
  title= {Generative adversarial nets},
  author={Goodfellow, I. and Pouget-Abadie, J. and Mirza, M. and Xu, B. and Warde-Farley, D. and Ozair, S. and Courville, A.  and Bengio, Y. },
  booktitle={Proceedings of the 27th Advances in Neural Information Processing Systems},
  pages={2672–2680},
  year={2014},
}

@article{2015Bow,
  author={Bowling, M. and Burch, N. and Johanson, M. and Tammelin, O. },
  journal={Science}, 
  title={Heads-up limit hold’em poker is solved}, 
  year={2015},
  volume={347},
  number={6218},
  pages={145–149},
  }

@article{Nedic2009Saddle,
  title={Subgradient methods for saddle-point problems},
  author={Nedi\'c, A. and Ozdaglar, A.},
  journal={Journal of optimization theory and applications},
  volume={142},
  number={1},
  pages={205-208},
  year={2009},
}

@article{2018Adrian,
  title={The Online Saddle Point Problem and Online Convex Optimization with Knapsacks},
  author={Cardoso, A. R. and Wang, H. and Xu, H. },
 journal={Mathematics of Operations Research},
 volume={50},
 number={1},
 pages={1-39},
  year={2025},
}

@article{1995Tseng,
  title={On linear convergence of iterative methods for the variational inequality problem},
  author={Tseng, P},
  journal={Journal of Computational and Applied Mathematics},
  year={1995},
  volume={60},
  number={1-2},
  pages={237–252},
}

@inproceedings{2020Mokhtari,
  title={A unified analysis of extra-gradient and optimistic gradient methods for saddle point problems: proximal point approach},
  booktitle={Proceedings of the 23rd International Conference on Artificial Intelligence and Statistics},
  author={Mokhtari, A. and Ozdaglar, A. and Pattathil, S},
  pages={1497–1507},
  year={2020},
}

@inproceedings{2019Liang,
  title={Interaction matters: a note on non-asymptotic local convergence of generative adversarial networks},
  booktitle={Proceedings of the 22nd International Conference on Artificial Intelligence and Statistics},
  author={ Liang, T  and  Stokes, J },
  pages={907–915},
  year={2019},
}

@article{Shalevshwartz2012Online,
  title={Online learning and online convex optimization},
  author={Shalev-Shwartz, S.},
  journal={Foundations and Trends in Machine Learning},
  volume={4},
  number={2},
  pages={107-194},
  year={2012},
}

@article{2010Distributed,
  title={Distributed stochastic subgradient projection algorithms for convex optimization},
  author={ Ram, S. S. and Nedi\'c, A. and Veeravalli, V.},
  journal={Journal of Optimization Theory $\&$ Applications},
  volume={147},
  number={3},
  pages={516-545},
  year={2010},
}

@inproceedings{2003Zinkevich,
  title={Online convex programming and generalized infinitesimal gradient ascent},
  author={ Zinkevich, M. },
  booktitle={Proceedings of the 20th International Conference on Machine Learning},
  pages={928-936},
  year={2003},
}

@article{2007Hazan,
  title={ Logarithmic regret algorithms for online convex optimization},
  author={E. Hazan and A. Agarwal and S. Kale},
  journal={Machine Learning},
  volume={69},
  number={2-3},
  pages={169-192},
  year={2007},
}

@book{2004Morgenstern,
title = {Theory of Games and Economic Behavior},
author = {Von Neumann, J. and Morgenstern, O.},
publisher = {Princeton University Press},
address = {Princeton},
year = {2004},
}

@inproceedings{2012projection,
  title={Projection-free online learning},
  author={Hazan, E. and Kale, S.},
  booktitle={Proceedings of the 29th International Conference on Machine Learning},
  pages={1843-1850},
  year={2012},
}

@article{2013Auto,
  title={Distributed convergence to {Nash} equilibria in two-network zero-sum games},
  author={B. Gharesifard and J. Cort\'es},
  journal={Automatica},
  volume={49},
  number={6},
  pages={1683–1692},
  year={2013},
}

@incollection{von2007theory,
  title={Theory of games and economic behavior: 60th anniversary commemorative edition},
  author={Von Neumann, John and Morgenstern, Oskar},
  booktitle={Theory of games and economic behavior},
  year={2007},
  publisher={Princeton university press},
}

@article{nash1950equilibrium,
  title={Equilibrium points in n-person games},
  author={Nash Jr, John F},
  journal={Proceedings of the national academy of sciences},
  volume={36},
  number={1},
  pages={48--49},
  year={1950},
  publisher={National Acad Sciences}
}

@article{daskalakis2018limit,
  title={The limit points of (optimistic) gradient descent in min-max optimization},
  author={Daskalakis, Constantinos and Panageas, Ioannis},
  journal={Advances in neural information processing systems},
  volume={31},
  year={2018}
}

@article{schafer2019competitive,
  title={Competitive gradient descent},
  author={Sch{\"a}fer, Florian and Anandkumar, Anima},
  journal={Advances in Neural Information Processing Systems},
  volume={32},
  year={2019}
}

@inproceedings{adolphs2019local,
  title={Local saddle point optimization: A curvature exploitation approach},
  author={Adolphs, Leonard and Daneshmand, Hadi and Lucchi, Aurelien and Hofmann, Thomas},
  booktitle={The 22nd International Conference on Artificial Intelligence and Statistics},
  pages={486--495},
  year={2019},
}

@article{hamedani2021primal,
  title={A primal-dual algorithm with line search for general convex-concave saddle point problems},
  author={Hamedani, Erfan Yazdandoost and Aybat, Necdet Serhat},
  journal={SIAM Journal on Optimization},
  volume={31},
  number={2},
  pages={1299--1329},
  year={2021},
}

@article{mokhtari2020convergence,
  title={Convergence rate of O(1/k) for optimistic gradient and extragradient methods in smooth convex-concave saddle point problems},
  author={Mokhtari, Aryan and Ozdaglar, Asuman E and Pattathil, Sarath},
  journal={SIAM Journal on Optimization},
  volume={30},
  number={4},
  pages={3230--3251},
  year={2020},
}

@inproceedings{yoon2021accelerated,
  title={Accelerated Algorithms for Smooth Convex-Concave Minimax Problems with O (1/k\^{} 2) Rate on Squared Gradient Norm},
  author={Yoon, TaeHo and Ryu, Ernest K},
  booktitle={International Conference on Machine Learning},
  pages={12098--12109},
  year={2021},
}

@inproceedings{liu2020min,
  title={Min-max optimization without gradients: Convergence and applications to black-box evasion and poisoning attacks},
  author={Liu, Sijia and Lu, Songtao and Chen, Xiangyi and Feng, Yao and Xu, Kaidi and Al-Dujaili, Abdullah and Hong, Mingyi and O’Reilly, Una-May},
  booktitle={International conference on machine learning},
  pages={6282--6293},
  year={2020},
}

@article{wang2023zeroth,
  title={Zeroth-order algorithms for nonconvex-strongly-concave minimax problems with improved complexities},
  author={Wang, Zhongruo and Balasubramanian, Krishnakumar and Ma, Shiqian and Razaviyayn, Meisam},
  journal={Journal of Global Optimization},
  volume={87},
  number={2},
  pages={709--740},
  year={2023},
}

@inproceedings{gidel2017frank,
  title={Frank-wolfe algorithms for saddle point problems},
  author={Gidel, Gauthier and Jebara, Tony and Lacoste-Julien, Simon},
  booktitle={Artificial Intelligence and Statistics},
  pages={362--371},
  year={2017},
}

@article{abernethy2017frank,
  title={On Frank-Wolfe and equilibrium computation},
  author={Abernethy, Jacob D and Wang, Jun-Kun},
  journal={Advances in Neural Information Processing Systems},
  volume={30},
  year={2017}
}

@inproceedings{abernethy2018faster,
  title={Faster rates for convex-concave games},
  author={Abernethy, Jacob and Lai, Kevin A and Levy, Kfir Y and Wang, Jun-Kun},
  booktitle={Conference On Learning Theory},
  pages={1595--1625},
  year={2018},
  organization={PMLR}
}

@inproceedings{hazan2017efficient,
  title={Efficient regret minimization in non-convex games},
  author={Hazan, Elad and Singh, Karan and Zhang, Cyril},
  booktitle={International Conference on Machine Learning},
  pages={1433--1441},
  year={2017},
}

@inproceedings{lin2020gradient,
  title={On gradient descent ascent for nonconvex-concave minimax problems},
  author={Lin, Tianyi and Jin, Chi and Jordan, Michael},
  booktitle={International Conference on Machine Learning},
  pages={6083--6093},
  year={2020},
  organization={PMLR}
}

@article{Kakutani1941,
  title={A generalization of brouwer’s fixed point theorem},
  author={Shizuo Kakutani},
  journal={Duke Mathematical Journal},
  volume={8},
  number={3},
  pages={457--459},
  year={1941},
}

@article{1928Neumann,
  title={Zur theorie der gesellschaftsspiele},
  author={Neumann, J},
  journal={Mathematische Annalen},
  volume={100},
  number={},
  pages={295--320},
  year={1928},
}

@inproceedings{2004Online,
  title={Online convex optimization in the bandit setting: gradient descent without a gradient},
  author={ Flaxman, A. D.  and  Kalai, A. T.  and  Mcmahan, H. B. },
  booktitle={Proceedings of the
16th Annual ACM-SIAM Symposium on Discrete Algorithms (SODA)},
  pages={385-394},
  year={2005},
}

@inproceedings{2018Dynamic,
  title={Dynamic Regret of Strongly Adaptive Methods},
  author={ Zhang, L.  and  Yang, T.  and  Zhou, Z.  },
  booktitle={Proceedings of the 35th International Conference on Machine Learning},
  pages={5882-5891},
  year={2018},
}

@inproceedings{2-xu2019online,
  title={An online saddle point optimization algorithm with regularization},
  author={Xu, Yue and Jiang, Ying and Xie, Xin and Li, Dequan},
  booktitle={IOP Conference Series: Materials Science and Engineering},
  volume={569},
  number={5},
  pages={052035},
  year={2019},
  organization={IOP Publishing}
}

@article{28-cardoso2019competing,
  title={Competing against equilibria in Zero-Sum Games with evolving payoffs},
  author={Cardoso, Adrian Rivera and Abernethy, Jacob and Wang, He and Xu, Huan},
  journal={arXiv preprint arXiv:1907.07723},
  year={2019}
}

@article{3-roy2019online,
  title={Online and bandit algorithms for nonstationary stochastic saddle-point optimization},
  author={Roy, Abhishek and Chen, Yifang and Balasubramanian, Krishnakumar and Mohapatra, Prasant},
  journal={arXiv preprint arXiv:1912.01698},
  year={2019}
}

@article{27-meng2023online,
  title={Online Saddle Point Problem and Online Convex-Concave Optimization},
  author={Meng, Qing-xin and Liu, Jian-wei},
  journal={arXiv preprint arXiv:2312.06957},
  year={2023},
}

@inproceedings{29-zhang2022no,
  title={No-regret learning in time-varying zero-sum games},
  author={Zhang, Mengxiao and Zhao, Peng and Luo, Haipeng and Zhou, Zhi-Hua},
  booktitle={Proceedings of the International Conference on Machine Learning},
  pages={26772--26808},
  year={2022},
}

@inproceedings{2021weilinear,
  title={Linear Last-iterate Convergence in Constrained Saddle-point Optimization},
  author={Wei, Chen-Yu and Lee, Chung-Wei and Zhang, Mengxiao and Luo, Haipeng},
  booktitle={Proceedings of the International Conference on Learning Representations},
  year={2021},
}

@INPROCEEDINGS{7963457,
  author={S. {Shahrampour} and A. {Jadbabaie}},
  booktitle={2017 American Control Conference},
  title={An online optimization approach for multi-agent tracking of dynamic parameters in the presence of adversarial noise},
  year={2017},
  volume={},
  number={},
  pages={3306-3311},}

@article{2019Introduction,
  title={Introduction to Online Convex Optimization},
  journal={Foundations and Trends in Optimization},
  volume={2},
  number={3-4},
  pages={157-325},
  author={ Hazan, E. },
  year={2016},
}

@article{liu2023forgetting,
  title={Forgetting-Factor Regrets for Online Convex Optimization},
  author={Liu, Yuhang and Zhao, Wenxiao and Yin, George},
  journal={IEEE Transactions on Automatic Control},
  volume={69},
  number={8},
  pages={5034-5048},
  year={2024},
}

@inproceedings{hazan2012online,
  title={Projection-free online learning},
  author={Hazan, E. and Kale, S.},
  booktitle={Proceedings of the 29th International Conference on Machine Learning},
  pages={1843-1850},
  year={2012},
}

\end{document}